\documentclass[12pt]{amsart}

\usepackage{amssymb}

\usepackage{xy}\xyoption{all}

\newtheorem{theorem}{Theorem}[section]

\newtheorem{lemma}[theorem]{Lemma}
\newtheorem{prop}[theorem]{Proposition}
\newtheorem{corollary}[theorem]{Corollary}
\newtheorem*{prop*}{Proposition}
\newtheorem*{lemma*}{Lemma}
\newtheorem*{claim*}{Claim}

\theoremstyle{definition}

\newtheorem{example}[theorem]{Example}
\newtheorem{remark}[theorem]{Remark}

\newtheorem{examples}[theorem]{Examples}

\newtheorem*{conjecture}{Conjectures}

\theoremstyle{remark}

\numberwithin{equation}{section}

\def\cstarlambda{C^\ast(\Lambda)}
\def\cstar{C^\ast}
\def\ZL0{{\Z^{\Lambda^0}}}

\def\Z{\mathbb{Z}}
\def\N{\mathbb{N}}
\def\T{\mathbb{T}}
\def\O{\mathcal{O}}
\def\F{\mathcal{F}}
\def\Q{\mathbb{Q}}

\newcommand{\Obj}{\operatorname{Obj}}
\newcommand{\coker}{\operatorname{coker}}
\newcommand{\img}{\operatorname{img}}
\newcommand{\Zmod}[1]{\Z/#1\Z}
\newcommand{\Circ}{\operatorname{Circ}}

\usepackage{pictexwd}

\newcommand{\sock}[4]{\beginpicture \setcoordinatesystem units <#1cm,#1cm>
\setplotarea x from 0 to 4, y from 0 to 4 \setlinear \plot 0 0 4 0 /
\plot 0 0 0 4 / \plot 0 2 4 2 / \plot 2 0 2 4 / \plot 0 4 2 4 /
\plot 4 2 4 0 / \put{#2} at 1 1 \put{#3}
 at 3 1 \put{#4} at 1 3 \endpicture}

\def\twodominoup{\beginpicture \setcoordinatesystem units <0.2cm,0.2cm>
\setplotarea x from 0 to 1, y from 0 to 2 \setlinear \plot 0 0 0 2 /
\plot 0 0 1 0 / \plot 1 0 1 2 / \plot 0 1 1 1 / \plot 0 2 1 2 /
\endpicture}

\def\threestair{\beginpicture \setcoordinatesystem units <0.2cm,0.2cm>
\setplotarea x from 0 to 3, y from 0 to 3 \setlinear \plot 0 0 3 0 /
\plot 0 1 3 1 / \plot 0 2 2 2 / \plot 0 3 1 3 / \plot 0 0 0 3 /
\plot 1 0 1 3 / \plot 2 0 2 2 / \plot 3 0 3 1 /
\endpicture}

\def\bootup{\beginpicture \setcoordinatesystem units <0.2cm,0.2cm>
\setplotarea x from 0 to 2, y from 0 to 3 \setlinear \plot 0 0 2 0 /
\plot 0 1 2 1 / \plot 0 2 1 2 / \plot 0 3 1 3 / \plot 0 0 0 3 /
\plot 1 0 1 3 / \plot 2 0 2 1 / \endpicture}

\begin{document}
\title[$2$-graphs from subshifts]{\boldmath{A family of $2$-graphs arising from\\ two-dimensional subshifts}}
\author[David Pask]{David Pask}
\address{David Pask, School  of Mathematics and Applied Statistics, University of Wollongong, NSW 2522, Australia}
\email{dpask@uow.edu.au}
\author[Iain Raeburn]{Iain~Raeburn}
\address{Iain Raeburn, School  of Mathematics and Applied Statistics, University of Wollongong, NSW 2522, Australia}
\email{raeburn@uow.edu.au}
\author[Natasha A. Weaver]{Natasha A. Weaver}
\address{Natasha A. Weaver, School of Mathematical and Physical Sciences,
University of Newcastle, NSW 2308, Australia}
\email{Natasha.Weaver@studentmail.newcastle.edu.au}

\thanks{This research was supported by the Australian Research Council, and Natasha Weaver was supported by an Australian Postgraduate Award. David Pask thanks Tom Ward and Klaus Schmidt for helpful discussions.}

\begin{abstract}
Higher-rank graphs (or $k$-graphs) were introduced by Kumjian and
Pask to provide combinatorial models for the higher-rank
Cuntz-Krieger $C^*$-algebras of Robertson and Steger. Here we
consider a family of finite $2$-graphs whose path spaces are
dynamical systems of algebraic origin, as studied by Schmidt and
others. We analyse the $C^*$-algebras of these $2$-graphs, find
criteria under which they are simple and purely infinite, and
compute their $K$-theory. We find examples whose $C^*$-algebras
satisfy the hypotheses of the classification theorem of Kirchberg
and Phillips, but are not isomorphic to the $C^*$-algebras of
ordinary directed graphs.
\end{abstract}

\maketitle \begin{center}\today \end{center}

\section{Introduction}

Higher-rank graphs (or $k$-graphs) were introduced by Kumjian and
Pask \cite{KP} to provide combinatorial models for the higher-rank
Cuntz-Krieger $C^*$-algebras of Robertson and Steger \cite{RS2}.
They have since  provided a fertile source of examples in
noncommutative geometry \cite{PRS1, PRS2, PZ}, and many important
operator algebras can be realised as the $C^*$-algebras of higher-rank graphs
\cite{KribsP, PRRS, DPY}. There has therefore been continuing
interest in finding new families of
$k$-graphs and analysing the structure of their $C^*$-algebras.

Every shift of finite type is equivalent to the backward shift
$\sigma$ on the two-sided infinite path space of a finite directed
graph \cite[Theorem~2.5]{LM}. The two-sided infinite path space
$\Lambda^\Delta$ of a finite $k$-graph $\Lambda$ introduced in
\cite{KP2} carries a set of $k$ commuting shifts $\sigma_i$, and
these are examples of the higher-dimensional shifts of finite type
studied by dynamicists. In this paper we consider a family of finite
$2$-graphs $\Lambda$ for which the path spaces
$(\Lambda^\Delta,\sigma_i)$ are dynamical systems of algebraic
origin, as studied by Schmidt and others \cite{Schmidt}. (A
particular motivating example for us was the system introduced by
Ledrappier in \cite{Ledrapp}.) We analyse the $C^*$-algebras of
these $2$-graphs, find criteria under which they are simple and
purely infinite, and compute their $K$-theory.

Each of our graphs $\Lambda$ is associated to a \emph{tile}, which
is a finite hereditary
subset $T$ of $\N^2$ containing the origin. We picture $T$ as a
collection of boxes into which we can put elements of the commutative ring $\Z/q\Z$, which we think of as an alphabet:
for example, we picture the \emph{sock} $T:=\{0,e_1,e_2\}$ as
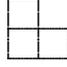
\[\sock{0.2}{}{}{}\]
The vertices in $\Lambda$ are copies of $T$ in which each box is
filled with elements of $\Z/q\Z$ which together satisfy a fixed equation in
$\Z/q\Z$; for example, the vertices in the \emph{Ledrappier graph} underlying
Ledrappier's system are copies of the sock filled with $0$s and $1$s
such that sum of the entries is $0 \pmod 2$. Paths in $\Lambda^*$
are diagrams covered by translates of $T$, filled in so that each
translate of $T$ is a valid vertex. Thus for example,
\begin{equation}\label{expath}
\beginpicture \setcoordinatesystem units <0.5cm,0.5cm>
\setplotarea x from 0 to 5, y from 0 to 4 \setlinear \plot 0 0 5 0 /
\plot 0 1 5 1 / \plot 0 2 5 2 / \plot 0 3 5 3 / \plot 0 4 4 4 /
\plot 0 0 0 4 / \plot 1 0 1 4 / \plot 2 0 2 4 / \plot 3 0 3 4 /
\plot 4 0 4 4 / \plot 5 0 5 3 /

\put{$0$} at 0.5 0.5 \put{$0$} at 1.5 0.5 \put{$1$} at 2.5 0.5
\put{$1$} at 3.5 0.5 \put{$0$} at 4.5 0.5 \put{$0$} at 0.5 1.5
\put{$1$} at 1.5 1.5 \put{$0$} at 2.5 1.5 \put{$1$} at 3.5 1.5
\put{$1$} at 4.5 1.5 \put{$1$} at 0.5 2.5 \put{$1$} at 1.5 2.5
\put{$1$} at 2.5 2.5 \put{$0$} at 3.5 2.5 \put{$1$} at 4.5 2.5
\put{$0$} at 0.5 3.5 \put{$0$} at 1.5 3.5 \put{$1$} at 2.5 3.5
\put{$1$} at 3.5 3.5
\endpicture
\end{equation}
represents a path $\lambda$ of degree $(3,2)$ in the Ledrappier
graph \[\text{from } s(\lambda)=\sock{0.2}{0}{1}{1} \text{ (the top RH one) }
\text{ to } r(\lambda)= \sock{0.2}{0}{0}{0} \text{ (the bottom LH one). }\] The infinite path
space $\Lambda^\Delta$ consists of similar diagrams covering the
entire plane, and the shifts $\sigma_1$ and $\sigma_2$ simply move
the diagram one row left and one column down respectively. So it is
easy to construct paths in these graphs, and thereby determine
properties of their $C^*$-algebras.

We begin with a short section in which we recall the basic
properties of $2$-graphs and their $C^*$-algebras. In
\S\ref{s_tiles}, we fix a set of ``basic data", which consists of a
tile $T$, an integer $q$ determining the alphabet, another integer
$t$ and a function $w:T\to \Z/q\Z$ which determines the equation relating the entries. We
describe the vertices and paths in our $2$-graph as functions from
$T$ and translates of $T$ into $\Z/q\Z$, and then we have to prove
that they form the objects and morphisms in a category satisfying
the axioms of a $2$-graph $\Lambda=\Lambda(T,q,t,w)$ (see
Theorem~\ref{thm_2graph}). In \S\ref{s_infpathsp}, we show that the
two-sided infinite path space of $\Lambda(T,q,t,w)$ is the
underlying space for a higher-dimensional shift of the sort studied
in dynamical systems (Theorem~\ref{thm_infpathsphomeo}).

In \S\ref{s_aperiodicity}, we show that, provided certain key values of the function $w$ are invertible, the $2$-graph $\Lambda(T,q,t,w)$ is aperiodic in the sense of Kumjian and Pask, so that the Cuntz-Krieger uniqueness theorem applies. We prove this using the recent formulation of Robertson and Sims \cite{RobS}, and as an intermediate step in the proof we show that $\Lambda(T,q,t,w)$ is always strongly connected in the sense that there are paths joining any two vertices. In \S\ref{s_cstaralgs}, we show that, under the same invertibility hypothesis on $w$, the $C^*$-algebra $C^*(\Lambda(T,q,0,w))$ is nuclear, simple and purely infinite. After our work in the previous section, this follows quickly from general results in \cite{RobS} and \cite{Sims} about the structure of $k$-graph algebras.
The main result, Theorem~{6.1}, implies that $C^*(\Lambda(T,q,0,w))$ is a Kirchberg algebra, and hence by the theorem of Kirchberg and Phillips is classifiable by its $K$-groups.

In \S\ref{s_kthry}, we compute the $K$-theory of $C^*(\Lambda(T,q,t,w))$, using the techniques developed by Robertson-Steger \cite{RSk} and Evans \cite{E}, which identify $K_0(C^*(\Lambda(T,q,t,w)))$ and $K_1(C^*(\Lambda(T,q,t,w)))$ in terms of the kernels and cokernels of certain integer matrices. Our $2$-graphs are finite but large, so we have used the computational algebra system \texttt{Magma} \cite{magma} to compute these kernels and cokernels. We have presented some of these results in Table~\ref{table_kthrytable}. These results have led us to make some general conjectures about the $K$-theory of our $2$-graphs, and in \S\ref{s_kthryresults} we prove two of these conjectures. Perhaps the most surprising result is that, under mild hypotheses, $K_0$ and $K_1$ have the same finite cardinality --- though our proof of this is indirect, and gives us no hint of whether the groups are actually isomorphic.

\subsection*{Visualisation convention}
We visualise a subset $S$ of $\N^2$ as the union of the unit squares
whose bottom left-hand corners belong to $S$, and a function $f:S\to
\Z$ as a diagram in which the number $f(i)$ is placed in the square
with bottom left-hand corner $i$. Thus the sock $T=\{0,e_1,e_2\}$ is
visualised as
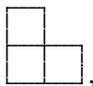
\[\sock{0.25}{}{}{}\,,\]
and the function $f:T\to \Z$ defined by $f(0)=3$, $f(e_1)=6$ and $f(e_2)=5$ as
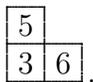
\[\sock{0.25}{3}{6}{5}\,.\]

\subsection*{Notation}
Let $\N=\{0,1,2,3,\ldots\}$ denote the monoid of natural numbers
under addition and let $\Z$ be the group of integers. For $k\geq 1$,
we view $\N^k$ as the set of morphisms in a
category with one object and composition
given by addition. We
write $n_i$ for the $i$th coordinate of $n\in\Z^k$, and $\{e_i\}$ for the usual basis of $\Z^k$. For $m,n\in\Z^k$
we say $m\leq n$ if $m_i\leq n_i$ for each $i$, and write $m\vee n$
and $m\wedge n$ for the coordinate-wise maximum and minimum.

\section{$2$-graphs}\label{s_2graphs}

Let $k$ be a positive integer. A
\emph{graph of rank $k$} or \emph{$k$-graph} is a pair $(\Lambda,d)$
consisting of a countable category $\Lambda$ and a functor
$d:\Lambda\to\mathbb{N}^k$, called the \emph{degree map}, satisfying
the \emph{factorisation property}: for every $\lambda\in\Lambda$ and
$m,n\in\N^k$ with $d(\lambda)=m+n$, there exist unique elements
$\mu,\nu\in\Lambda$ such that $d(\mu)=m$, $d(\nu)=n$ and
$\lambda=\mu\nu$. In practice, we drop the degree map from the notation.

We refer to the morphisms in $\Lambda$ as \emph{paths} and the
objects as \emph{vertices}. If $\lambda\in \Lambda$ satisfies
$d(\lambda)=n$ we say $\lambda$ has \emph{degree} $n$; we write
$\Lambda^n:=\{\lambda\in\Lambda:d(\lambda)=n\}$. All $k$-graphs in
this paper are \emph{finite} in the sense that each $\Lambda^n$ is a
finite set, and have no sources in the sense that for each vertex $v$ and each $n\in \N^k$, there is at least one $\lambda\in \Lambda^n$ with $r(\lambda)=v$.

The factorisation property has several consequences. First, it
implies that for each vertex $v$ the identity morphism $\iota_v$ is
the only morphism of degree $0$ from $v$ to $v$, so that we can
identify $\Obj(\Lambda)$ with $\Lambda^0$. We then write $v\Lambda$,
for example, to mean the set of paths $\lambda$ with $r(\lambda)=v$.
Second, it implies that for every triple $m,n,p\in \N^k$ satisfying
$0\leq m\leq n\leq p$ and $\lambda\in\Lambda^p$ with $d(\lambda)=p$,
there are unique \emph{segments} $\lambda(0,m)\in\Lambda^m$,
$\lambda(m,n)\in\Lambda^{n-m}$, $\lambda(n,p)\in\Lambda^{p-n}$ such
that $\lambda=\lambda(0,m)\lambda(m,n)\lambda(n,p)$. The paths
$\lambda(m,m)$ have degree $0$, and hence are vertices; in the
literature it is common to write $\lambda(m):=\lambda(m,m)$, but we
will refrain from doing this as $\lambda(m)$ will have another more
natural meaning.

\begin{figure}
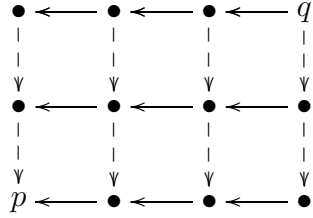

\centerline{
\xygraph{
{q}="v43":[l]{\bullet}="v33":[l]{\bullet}="v23":[l]{\bullet}="v13"
"v43":@{-->}[d]{\bullet}="v42":@{-->}[d]{\bullet}="v41"
"v33":@{-->}[d]{\bullet}="v32":@{-->}[d]{\bullet}="v31"
"v23":@{-->}[d]{\bullet}="v22":@{-->}[d]{\bullet}="v21"
"v13":@{-->}[d]{\bullet}="v12":@{-->}[d]{p}="v11"
"v43":"v33":"v23":"v13"
"v42":"v32":"v22":"v12"
"v41":"v31":"v21":"v11"
}}
\caption{A path of degree (3,2).}\label{figpath}
\end{figure}

\begin{figure}
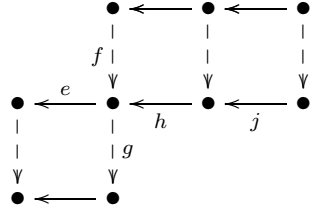
\label{figcompose}
\centerline{ \xygraph{
{\bullet}="v43":[l]{\bullet}="v33":[l]{\bullet}="v23":@{}[l]{}="v13"
"v43":@{-->}[d]{\bullet}="v42" "v33":@{-->}[d]{\bullet}="v32"
"v23":@{-->}[d]{\bullet}="v22"_f:@{-->}[d]{\bullet}="v21"^g
"v13":@{}[d]{\bullet}="v12":@{-->}[d]{\bullet}="v11"
"v43":"v33":"v23" "v42":"v32"^j:"v22"^h:"v12"_e "v21":"v11" }}
\caption{Composing paths.}
\end{figure}

In this paper we are primarily interested in $2$-graphs, so $k$ is
usually $2$. We visualise a $2$-graph as a directed bicoloured graph
with vertex set $\Lambda^0$ in which the elements $\beta$ of
$\Lambda^{e_1}$ are represented by blue edges from $s(\beta)\in
\Lambda^0$ to $r(\beta)\in \Lambda^0$, and elements of
$\Lambda^{e_2}$ as red edges. (In print we use black curves to
represent blue edges and dashed curves for red edges.) This
bicoloured graph is called the \emph{skeleton} of $\Lambda$.
Applying the factorisation property to $(1,1)=e_1+e_2=e_2+e_1$ gives
a bijection between the blue-red paths of length $2$ and the
red-blue paths of length $2$. We then visualise a path of degree
$(1,1)$ as a square
\begin{equation}
\label{exsquare}
\xygraph{{\bullet}="v11":@{-->}[d]{\bullet}="v10"^h:[l]{\bullet}="v00"^g
"v11":[l]{\bullet}="v01"_f:@{-->}"v00"_e}
\end{equation}
in which the bijection matches up the blue-red path $gh$ with the
red-blue path $ef$, so that $gh=ef$ are the two factorisations of
the path of degree $(1,1)$. It turns out (though we shall not rely
on this fact in this paper) that a $2$-graph is completely
determined by a collection $C$ of squares \eqref{exsquare} in which
each blue-red and each red-blue path occur exactly once. The paths
of degree $(3,2)$ from $w$ to $v$, for example, then consist of
copies of the rectangle in Figure~\ref{figpath}
pasted round the blue-red graph, so that $q$ lands on $w$, $p$ lands
on $v$, and each constituent square is one of the given collection
$C$. Composition of paths involves taking the convex hull: if
$d(\lambda)=(1,1)$ and $d(\mu)=(1,2)$, for example, then
$\lambda\mu$ is obtained by filling in the corners of the diagram in
Figure~\ref{figcompose}
with squares from $C$, which can be done in exactly one way (there
is only one square fitting $ef$, for example, so this has to be the square in the top left corner).

If $\Lambda$ is a finite $k$-graph with no sources, then a
\emph{Cuntz-Krieger $\Lambda$-family} is a collection of partial
isometries $\{S_\lambda:\lambda\in \Lambda\}$ (either operators in a
Hilbert space or elements of an abstract $C^*$-algebra) satisfying
the following \emph{Cuntz-Krieger relations}:
\begin{enumerate}
\item $\{S_v:v\in\Lambda^0\}$ is a family of mutually orthogonal
projections,
\smallskip
\item $S_{\lambda\mu}=S_\lambda S_\mu$ for all
$\lambda,\mu\in \Lambda$ such that $s(\lambda)=r(\mu)$,
\smallskip
 \item
$S_\lambda^\ast S_\lambda=S_{s(\lambda)}$ for all $\lambda\in
\Lambda$, and
\smallskip \item $S_v=\sum_{\lambda\in v\Lambda^n}S_\lambda S_\lambda^\ast$ for all $v\in \Lambda^0$ and $n\in \N^k$.
\end{enumerate}
The $C^*$-algebra of $\Lambda$ is the $C^*$-algebra
$C^\ast(\Lambda)$ generated by a universal Cuntz-Krieger
$\Lambda$-family $\{s_\lambda:\lambda\in \Lambda\}$. The basic facts
about these $C^*$-algebras are discussed in \cite{KP}, \cite{RSY03}
and \cite[Chapter~10]{R}, for example.

\section{Tiles and $2$-graphs}\label{s_tiles}

There are four variables in our main construction of $2$-graphs.
Recall that a subset $T$ of $\N^2$ is \emph{hereditary} if $j\in T$
and $0\leq i\leq j$ imply $i\in T$. The variables are:
\begin{itemize}
\item a \emph{tile} $T$, which is a hereditary subset
of $\N^2$ with finite cardinality $|T|$;

\item an \emph{alphabet} $\{0,1,\ldots,q-1\}$, where $q\geq 2$ is an integer; we view the alphabet as a commutative ring by identifying it with $\Zmod{q}$ in the obvious way;

\item an element $t$ of the alphabet, called the \emph{trace}; and

\item a function $w:T\to \{0,1,\ldots,q-1\}$ called the \emph{rule}.
\end{itemize}
For the rest of the section, we fix the \emph{basic data} $(T,q,t,w)$.

\begin{example}\label{socktile}
For the $2$-graph underlying the Ledrappier system, the basic data
consists of the sock tile $T=\{0,e_1,e_2\}$, $q=2$, $t=0$ and the
constant function $w\equiv 1$.
\end{example}

The vertex set in our $2$-graph will be
\[\Lambda^0 =\Big\{v:T\to\Z/q\Z : \sum_{i\in T}w(i)v(i)=t\pmod{q}\Big\}.\]
The Ledrappier graph, for example, has four vertices $a,b,c,d$ visualised as\begin{equation}\label{Ledrapvert}
\sock{0.25}{0}{0}{0} \hspace{1cm}\sock{0.25}{0}{1}{1} \hspace{1cm}\sock{0.25}{1}{0}{1} \hspace{1cm} \sock{0.25}{1}{1}{0}\,.
\end{equation}

To describe the paths, we need some notation. Let
$(c_1,c_2):=\bigvee\{i:i\in T\}$, so that the longest row in $T$
(the bottom one) has $c_1+1$ boxes and the highest column in $T$
(the left-hand one) has $c_2+1$ boxes. For $S\subset \Z^2$ and $n\in
\Z^2$, we let $S+n=\{i+n:i\in S\}$ denote the translate of $S$ by
$n$, and we set $T(n):=\bigcup_{0\leq m\leq n}T+m$. When we
visualise $T(n)$ using our convention, it looks like a
$(c_1+1+n_1)\times (c_2+1+n_2)$ rectangle of boxes with a bite taken
out of the top right-hand corner. If $f:S\to\Z/q\Z$ is a function
defined on a subset $S$ of $\N^2$ containing $T+n$, then we define
$f|_{T+n}:T\to\Z/q\Z$ by
\begin{equation}
\label{defrestriction}f|_{T+n} (i) = f(i+n) \text{ for } i\in T.
\end{equation}
A \emph{path of degree $n$} is a function $\lambda:T(n)\to \Z/q\Z$
such that $\lambda|_{T+m}$ is a vertex for $0\leq m\leq n$; then
$\lambda$ has \emph{source} $s(\lambda)=\lambda|_{T+n}$ and
\emph{range} $r(\lambda)=\lambda|_T$. Thus, for example, the diagram
\eqref{expath} is the visualisation of a path of degree $(3,2)$ in
the Ledrappier graph based on the sock tile.

Notice that the function $f|_{T+n}$ defined in
\eqref{defrestriction} is not a simple restriction: because our
tiles all have their bottom left-hand corner at the origin, we need
to translate by $n$ on the right-hand side. We need to use a similar
convention when we define the segments appearing in the
factorisations of paths. For $\lambda\in\Lambda^p$ and $0\leq m\leq
n\leq p$, the \emph{segment} $\lambda(m,n)$ is the path of degree
$n-m$ defined by
\[ \lambda(m,n)(i) = \lambda(m+i) \text{ for } i\in T(n-m).\]
In particular, $\lambda(m,m)$ is the vertex $\lambda|_{T+m}$.

We want $\Lambda^\ast:=\bigcup_{n\in \N^2}\Lambda^n$ to be the
morphisms in a category, and so we have to define composition. To
make this work, we need to make an assumption on the rule $w$.  We
say that the rule $w$ \emph{has invertible corners} if $w(c_1e_1)$
and $w(c_2e_2)$ are invertible elements of the ring $\Zmod{q}$. The
next proposition tells us that there is exactly one candidate for the
composition of two paths.

\begin{prop}\label{prop_pathcomposition}
Suppose we have basic data $(T,q,t,w)$ and the rule $w$ has invertible corners. Suppose
$\mu\in\Lambda^m$ and $\nu\in\Lambda^n$ satisfy $s(\mu)=r(\nu)$. Then there is a
unique path $\lambda\in \Lambda^{m+n}$ such that
\begin{equation} \lambda(0,m)=\mu\ \text{ and
}\ \lambda(m,m+n)=\nu. \label{eqn_pathcompcondn}
\end{equation}
\end{prop}

Notice that Equation~\eqref{eqn_pathcompcondn} defines $\lambda$
uniquely on $T(m)\cup (T(n)+m)$, so our problem is to show that
there is a unique function $\lambda':T(m+n)\to \Zmod{q}$ such that
$\lambda'|_{T(m)\cup (T(n)+m)}=\lambda$ and $\lambda'|_{T+k}$
belongs to $\Lambda^0$ for every $k$ such that $0\leq k\leq m+n$;
since $\mu$ and $\nu$ are paths and $\lambda'$ extends $\lambda$, we
already know this for $k$ such that $T+k\subset T(m)\cup (T(n)+m)$.

Our strategy is to extend $\lambda$ from $T(m)\cup (T(n)+m)$ to
$T(m+n)$ by adding one point at a time in such a way that there is
only one possible value for $\lambda$ at the new point. The next
lemma tells us how to do this.  It depends crucially on the
assumption that the rule has invertible corners (see
Remark~\ref{whyinvcorners}), and it fails spectacularly for
tiles in $\N^k$ when $k>2$ (see Remark~\ref{whynot3d}).

\begin{lemma}\label{oneatatime}
Suppose that $l\in \N^2$ and $S$ is a subset of $\N^2$ containing
$T+l-e_2$ and $T+l+e_1$, and $\lambda:S\to \Zmod{q}$ is a function
such that $\lambda|_{T+k}$ belongs to $\Lambda^0$ for every $k\in
\N^2$ such that $T+k\subset S$. Then there is a unique function
\[
\lambda':S':=S\cup\{l+e_1-e_2+c_1e_1\}\to\Zmod{q}
\]
such that $\lambda'|_S=\lambda$ and $\lambda'|_{T+k}$ belongs to
$\Lambda^0$ for every $k\in \N^2$ such that $T+k\subset S'$.
\end{lemma}

\begin{proof}
If $l+e_1-e_2+c_1e_1$ belongs to $S$, there is nothing to do. So we
suppose $l+e_1-e_2+c_1e_1$ is not in $S$. Let $i\in T\backslash
\{c_1e_1\}$. Then either $i_2=0$ and $i_1<c_1$, in which case
$i+e_1$ belongs to $T$ and $l+e_1-e_2+i$ belongs to $T+l-e_2$, or
$i_2>0$, in which case $i-e_2\in T$ and $l+e_1-e_2+i$ belongs to
$T+l+e_1$. So $(T\backslash \{c_1e_1\})+l+e_1-e_2$ is contained in
the domain $S$ of $\lambda$, and $l+e_1-e_2+c_1e_1$ is the only
point of $T+l+e_1-e_2$ which is not in $S$. Thus we can define
$\lambda'|_S=\lambda$ and
\begin{equation}\label{deflambda'}
\lambda'(l+e_1-e_2+c_1e_1):= w(c_1e_1)^{-1}\Big(t-\sum_{i\in
T\backslash \{c_1e_1\}}w(i)\lambda(l+e_1-e_2+i)\Big).
\end{equation}
If $T+k\subset S'$, then either $T+k\subset S$, in which case
$\lambda'|_{T+k}=\lambda|_{T+k}\in \Lambda^0$, or $k=l+e_1-e_2$, in
which case \eqref{deflambda'} implies $\lambda'|_{T+k}\in
\Lambda^0$. No other value of $\lambda'(l+e_1-e_2+c_1e_1)$ would
give $\lambda'|_{T+l+e_1-e_2}\in \Lambda^0$, so this function $\lambda'$ is
the only one with the required property.
\end{proof}

\begin{proof}[Proof of Proposition~\ref{prop_pathcomposition}]
The region $T(m+n)$ is obtained from $T(m)\cup (T(n)+m)$ by adding two rectangles
\begin{align*}
BR&:=\{j\in \N^2:c_1+m_1<j_1\leq c_1+m_1+n_1,\ 0\leq j_2<m_2\},\text{ and}\\
UL&:=\{j\in \N^2:0\leq j_1<m_1,\ c_2+m_2<j_2\leq c_2+m_2+n_2\}.
\end{align*}
We order the bottom right rectangle $BR$ lexicographically, first
going down the column $j_1=c_1+m_1+1$, then down the column
$j_1=c_1+m_1+2$, and so on. We then apply Lemma~\ref{oneatatime} to
each $j$ in order: when we come to define $\lambda(j)$, we have
already defined $\lambda(i)$ for every $i$ above and to the left of
$j$, and with $l:=j-e_1+e_2-c_1e_1$, $\lambda$ is defined on both
$T+l-e_2$ and $T+l+e_1$. Since there is only one possible value of
$\lambda(j)$ at each stage, there is only one way to extend
$\lambda$ to $T(m)\cup (T(n)+m)\cup BR$.

To see that $\lambda$ extends uniquely to $UL$, we can either run
the mirror image of this argument in the rectangle $UL$, or reflect
everything in the line $n_1=n_2$ and apply what we have just proved.
\end{proof}

\begin{theorem}\label{thm_2graph}
Suppose we have basic data $(T,q,t,w)$ and the rule $w$ has
invertible corners. Say that $\mu\in\Lambda^m$ and $\nu\in\Lambda^n$
are \emph{composable} if $s(\mu)=r(\nu)$, and define the
\emph{composition} $\mu\nu$ to be the unique path $\lambda$
satisfying \eqref{eqn_pathcompcondn}. Define $d:\Lambda\to\N^2$ by
$d(\lambda)=n$ for $\lambda\in\Lambda^n$. Then, with $\Lambda^0$,
$\Lambda^*$, $r$ and $s$ defined at the beginning of the section,
$\Lambda(T,q,t,w):=((\Lambda^0,\Lambda^*,r,s),d)$ is a $2$-graph.
\end{theorem}

\begin{proof}
We can view a vertex $v\in\Lambda^0$ as a path of degree $0$; then
$\lambda$ has the property \eqref{eqn_pathcompcondn} which
characterises $r(\lambda)\lambda$ and $\lambda s(\lambda)$, so $v$
has the properties required of the identity morphism at $v$. For
$\mu\in\Lambda^m$, $\nu \in\Lambda^n$ with $s(\mu)=r(\nu)$,
\eqref{eqn_pathcompcondn} implies that
$r(\mu\nu)=(\mu\nu)|_T=\mu|_T=r(\mu)$ and
\[
s(\mu\nu)=(\mu\nu)|_{T+m+n}=(\mu\nu)(m,m+n)|_{T+n}=\nu|_{T+n}=s(\nu).
\]
To prove that $\Lambda$ is a category, it remains to show that composition is associative.

Suppose $\mu\in\Lambda^m$, $\nu\in\Lambda^n$ and $\rho\in\Lambda^p$ satisfy
$s(\mu)=r(\nu)$ and $s(\nu)=r(\rho)$. For $i\in T(n)$, we have
\begin{align*}
((\mu\nu)\rho)(m,m+n+p)(i)&= ((\mu\nu)\rho)(i+m)=((\mu\nu)\rho)(0,m+n)(i+m) \\
&= (\mu\nu)(i+m)= (\mu\nu)(m,m+n)(i)=\nu(i),
\end{align*}
and for $i\in T(p)$ we have 
\[
((\mu\nu)\rho)(m,m+n+p)(i+n)=((\mu\nu)\rho)(i+m+n)=((\mu\nu)\rho)(m+n,m+n+p)(i)=\rho(i).
\]
Thus $((\mu\nu)\rho)(m,m+n+p)(0,n)=\nu$ and
$((\mu\nu)\rho)(m,m+n+p)(n,n+p)=\rho$, and hence
$((\mu\nu)\rho)(m,m+n+p)=\nu\rho$. On the other hand, for $i\in
T(m)$, we have
\begin{align*}
((\mu\nu)\rho)(0,m)(i)&= ((\mu\nu)\rho)(i)=((\mu\nu)\rho)(0,m+n)(i) \\
&= (\mu\nu)(i)= (\mu\nu)(0,m)(i)= \mu(i),
\end{align*}
so $((\mu\nu)\rho)(0,m)=\mu$. Thus $(\mu\nu)\rho$ has the property
which characterises $\mu(\nu\rho)$, and we have
$(\mu\nu)\rho=\mu(\nu\rho)$.

We have now shown that $\Lambda$ is a category, and it is countable
because each $\Lambda^n$ is finite. The map $d:\Lambda\to\N^2$
satisfies $d(\mu\nu)=d(\mu)+d(\nu)$ and hence is a functor, so it
remains to verify that $d$ has the factorisation property. But this
is easy: given $\lambda\in\Lambda^{m+n}$, the paths
$\mu:=\lambda(0,m)$ and $\nu:=\lambda(m,m+n)$ are the only ones
which can satisfy $\lambda=\mu\nu$.
\end{proof}

To visualise the $2$-graph $\Lambda(T,q,t,w)$, we draw its skeleton.
This skeleton has a few special properties.

\begin{prop}\label{prop_skeleton}
Suppose we have basic data $(T,q,t,w)$ and the rule $w$ has
invertible corners. Then $\Lambda=\Lambda(T,q,t,w)$ satisfies
\begin{itemize}
\item[\textnormal{(a)}] $|\Lambda^0|=q^{|T|-1}$;
\smallskip
\item[\textnormal{(b)}] for $v,u\in \Lambda^0$, $v\Lambda^{e_i}u$ is non-empty if and only if
\begin{equation}\label{existedge}
v(m)=u(m-e_i)\ \text{ for every $m\in T\cap (T+e_i)$,}
\end{equation}
in which case $|v\Lambda^{e_i}u|=1$;
\smallskip
\item[\textnormal{(c)}] $|v\Lambda^{e_1}|=|\Lambda^{e_1}v|=q^{c_2}$ and $|v\Lambda^{e_2}|=|\Lambda^{e_2}v|=q^{c_1}$ for every $v\in \Lambda^0$.
\end{itemize}
\end{prop}

\begin{proof}
There are $q^{|T|-1}$ functions $v:T\backslash\{c_1e_1\}\to \Z/q\Z$,
and each defines a unique vertex $v$ by setting
\[
v(c_1e_1)=w(c_1e_1)^{-1}\Big(t-\sum_{i\in T\backslash
\{c_1e_1\}}w(i)v(i)\Big).
\]
This gives (a). For (b), note that if $\beta\in v\Lambda^{e_i}u$ and $m\in T\cap (T+e_i)$, then
\[
v(m)=\beta|_T(m)=\beta|_{T+e_i}(m-e_i)=u(m-e_i).
\]
Conversely, if $u,v$ satisfy \eqref{existedge}, then we can define $\beta:T(e_i)\to \Z/q\Z$ by
\[
\beta(m)=\begin{cases}v(m)&\text{for $m\in T$,}\\
u(m-e_i)&\text{for $m\in (T+e_i)\backslash T$,}\end{cases}
\]
and \eqref{existedge} says that $\beta|_{T+e_i}=u$. The constraints
$\beta|_T=v$ and $\beta|_{T+e_i}=u$ completely determine $\beta$, so
$|v\Lambda^{e_i}u|=1$.

To see (c), note that an edge $\beta\in\Lambda^{e_1}v$ has
$\beta|_{T+e_1}$ determined by $v$. The remainder
$T(e_1)\setminus(T+e_1)$ is the first column of $T(e_1)$, which has
$c_2+1$ entries. There are $q^{c_2}$ ways of filling in the bottom
$c_2$ squares, and then the top entry is determined by
\[
\beta(c_2e_2)=w(c_2e_2)^{-1}\Big(t-\sum_{i\in T\backslash
\{c_2e_2\}}w(i)\beta(i)\Big).
\]
Thus $|\Lambda^{e_1}v|=q^{c_2}$. On the other hand, an edge
$\beta\in v\Lambda^{e_1}$ has $\beta|_T=v$, and $T(e_1)\backslash T$
also has $c_2+1$ squares. We can fill in all the squares except
$c_1e_1+e_1$ arbitrarily in $q^{c_2}$ ways, and then
$\beta(c_1e_1+e_1)$ is determined by
\[
\beta(c_1e_1+e_1)=w(c_1e_1)^{-1}\Big(t-\sum_{i\in T\backslash
\{c_1e_1\}}w(i)\beta(i+e_1)\Big).
\]
The facts about the red edges follow by symmetry.
\end{proof}

\begin{figure}
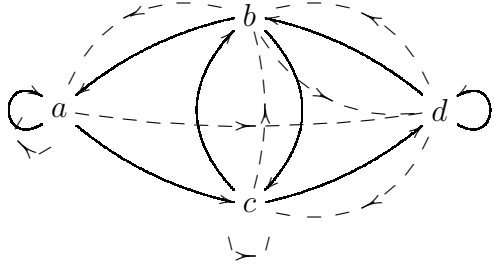
\label{Ledrapskeleton}
\centerline{\xygraph{{a}="a":@/_/[drr]c="d"_{}:@/_/[urr]
d="b"^{}:@/_/[ull]b="c"_{} "c":@/_/"a"^{}
"c":@/^{20pt}/"d"^{} "d":@/^{20pt}/"c"^{} "a":@(dl,ul)"a"^{}
"b":@(dr,ur)"b"^{} "c":@{--}|@{>}@/_{20pt}/"a"^{}
"b":@{--}|@{>}@/_{20pt}/"c"_{} "c":@{--}|@{>}@/_{15pt}/"b"^{}
"a":@{--}|@{>}@(d,l)"a"^{} "d":@{--}|@{>}@(dl,dr)"d"_{}
"d":@{--}|@{>}@/_/"c"_{} "a":@{--}|@{>}@/_/"b"^{}
"b":@{--}|@{>}@/^{20pt}/"d"_{}}}
\caption{The skeleton of the Ledrappier graph.}
\end{figure}

\begin{example}
The \emph{Ledrappier graph} $L(\sock{0.06}{}{}{})$ is the $2$-graph
constructed from the basic data consisting of the sock tile $T$,
$q=2$, $t=0$ and $w\equiv 1$. It has four vertices $a,b,c,d$ listed
in \eqref{Ledrapvert}. Examples of a blue edge (with range $b$ and
source $d$) and a red edge (with range $a$ and source $b$) are
visualised by
\[\beginpicture
\setcoordinatesystem units <0.25cm,0.25cm> \setplotarea x from 0 to
9, y from -1 to 4 \setlinear \plot 0 0 4 0 / \plot 0 0 0 4 / \plot 0
2 4 2 / \plot 2 0 2 4 / \plot 0 4 2 4 / \plot 4 2 4 0 / \put{$$} at
2 -1 \put{$0$} at 1 1 \put{$1$} at 1 3 \put{$1$} at 3 1
\plot 2 0 6 0 / \plot 2 0 2 4 /
\plot 2 2 6 2 / \plot 4 0 4 4 / \plot 2 4 4 4 / \plot
6 2 6 0 / \put{$$} at 4 -1 \put{$0$} at 3 3 \put{$1$} at 5 1 \put{and} at 9 3
\endpicture
\hspace{0.2cm}
\beginpicture 
\setcoordinatesystem units <0.25cm,0.25cm> \setplotarea x from 0 to
5, y from 0 to 6 \setlinear \plot 0 0 4 0 / \plot 0 0 0 4 / \plot 0
2 4 2 / \plot 2 0 2 4 / \plot 0 4 2 4 / \plot 4 2 4 0 / \put{$$} at
-1 2 \put{$0$} at 1 1 \put{$0$} at 1 3 \put{$0$} at 3 1
\plot 0 2 4 2 / \plot 0 2 0 6 /
\plot 0 4 4 4 / \plot 2 2 2 6 / \plot 0 6 2 6 / \plot 4
4 4 2 / \put{$$} at -1 4 \put{$1$} at 3 3 \put{$1$} at 1 5
\endpicture.\]
The skeleton of $L(\sock{0.06}{}{}{})$ is the $2$-coloured graph in Figure~\ref{Ledrapskeleton}.
\end{example}

\begin{remark}\label{whyinvcorners}
If we start with a rule which does not have invertible corners, then
we can still draw a bicoloured graph, which may or
may not be the skeleton of a $2$-graph. For example, suppose $T$ is the
sock, $w(0)=w(e_1)=1$ and $w(e_2)=0$. Then there is
exactly one blue-red path  between each pair of vertices, but there
are sometimes two and sometimes no red-blue paths, so the bicoloured graph cannot be
the skeleton of a $2$-graph. On the other hand, if we use the zero
rule $w\equiv 0$ on the sock, then there are two blue-red paths and
two red-blue paths between each pair of vertices, so there are
bijections between the sets of blue-red and red-blue paths, each
of which determines a potentially different $2$-graph. This is
reminiscent of the $2$-graphs $\F_\theta^+$ studied in \cite{DPY},
which have a single vertex, and are thus completely determined by a
permutation $\theta$ of a finite set
$\{1,\cdots,m\}\times\{1,\cdots,n\}$.

We will observe in Remark~\ref{whyinvert2} that when we have to make
choices to define a factorisation property, the correspondence
between $2$-graphs and shifts breaks down.
\end{remark}

\begin{remark}\label{whynot3d}
When we start with a tile $T$ which is a finite hereditary subset of
$\N^3$, we can construct a tricoloured graph, but this will not
completely determine a $3$-graph because Lemma~\ref{oneatatime}
fails. The crux of the proof of Lemma~\ref{oneatatime} (when
$l=e_2$) is that the set $T(e_1+e_2)\backslash (T(e_2)\cup
(T(e_1)+e_2))$, consists of the single point $(c_1+1)e_1$. When we
consider the tile $T=\{0,e_1,e_2,e_3\}$, which is a natural
$3$-dimensional analogue of the sock, we have
\[
T(e_1+e_2)\backslash (T(e_2)\cup (T(e_1)+e_2))=\{2e_1, e_1+e_3\}.
\]
If we use a single rule $w$ with invertible corners to define our
vertices and paths, then there is still more than one way to fill in
the two empty cubes. So one would have to impose more than one rule
to get a uniquely defined red-blue factorisation of a blue-red path.
However, the number of empty cubes to be filled depends on the
dimensions of the original tile, so just one extra rule is not
enough in general.
\end{remark}

\section{Connections with shift spaces}\label{s_infpathsp}

To make contact with the dynamics literature, we consider basic data
$(T,q,0,w)$. We denote by $R^q_2 = \Zmod{q}[u_1^{\pm 1},u_2^{\pm
1}]$ the commutative ring of Laurent polynomials in $u_1$, $u_2$
over the ring $\Zmod{q}$, and define $g=g_{T,w}\in R^q_2$ by
\[
g_{T,w}=\sum_{m\in T} w(m)u^m.
\]
The shift space $\Omega=\Omega^{R^q_2/(g)}$ is defined in \cite[page~719]{KitS4} as
\begin{align}
\Omega&= \Big\{ f=(f(n))\in (\Zmod{q})^{{\Z}^2}:\sum_{i\in T} w(i)f(i+n)=0 \pmod{q} \text{ for } n\in \Z^2 \Big\}\notag \\
&= \Big\{ f:\Z^2 \to\Zmod{q}:\sum_{i\in T} w(i)f|_{T+n}(i)=
0\pmod{q} \text{ for } n\in\Z^2 \Big\}\label{defX}.
\end{align}
This is a compact subspace of $( \Zmod{q} )^{\Z^2}$ in the product
topology, and carries an action of $\Z^2$ defined by
$({\alpha_p}f)(n)=f(n+p)$.

The two-sided path space of a $k$-graph $\Lambda$ was introduced in
\cite[\S 3]{KP2}; here we consider the case where $\Lambda$ is a
finite $2$-graph. Let
\[
\Delta=\{ (m,n):m,n\in\Z^2, m\le n \};
\]
with $r(m,n)=m$, $s(m,n)=n$ and $d(m,n)=n-m$, $(\Delta,d)$ becomes a
$2$-graph. The two-sided infinite path space of $\Lambda$ is
\[
\Lambda^\Delta=\{ x:\Delta\to\Lambda : x \text{ is a
degree-preserving functor} \}.
\]
It is shown in \cite{KP2} that $\Lambda^\Delta$ has a locally
compact (metric) topology with basic open sets
\begin{equation*} \label{tstopdef}
Z(\lambda,n)= \{ x\in \Lambda^\Delta : x(n,n+d(\lambda))=\lambda \},
\end{equation*}
for $\lambda\in\Lambda$ and $n\in\Z^2$. Since we are assuming that $\Lambda$ is finite, $\Lambda^\Delta$ is compact.

Now for each $p \in \Z^2$ we define $\sigma_p : \Lambda^\Delta \to \Lambda^\Delta$  by
\[ \sigma_p (x)(m,n) = x(m+p,n+p). \]
Observe that for all $n,p\in\Z^2$ and $\lambda\in\Lambda$ we have
$\sigma_p(Z(\lambda,n))=Z(\lambda,n+p)$, so $\sigma_p$ is a
homeomorphism for every $p\in\Z^2$, and $\sigma$ is an action of
$\Z^2$ on $\Lambda^\Delta$.

\begin{theorem}\label{thm_infpathsphomeo}
Suppose we have basic data $(T,q,0,w)$ and $w$ has invertible
corners, let $\Lambda:=\Lambda(T,q,0,w)$ be the associated
$2$-graph, and define $\Omega$ as in \eqref{defX}. Then there is a
homeomorphism $h:\Lambda^\Delta\to \Omega$ such that $\alpha_p\circ
h=h\circ\sigma_p$.
\end{theorem}

\begin{proof}
Define $h:\Lambda^\Delta\to(\Z/q\Z)^{\Z^2}$ by
\begin{equation*}
h(x)(i)=x(i,i)(0) \text{ for } x\in\Lambda^\Delta,i\in\Z^2.
\end{equation*}
Let $j\in T$. Since $x(i-j,i)$ is a path in $\Lambda$, it is a well-defined function on $T(j)$ and
\begin{equation}\label{eqn_red}
x(i,i)(0)=x(i-j,i)|_{T+j}(0)=x(i-j,i)(j+0)=x(i-j,i-j)(j).
\end{equation}
Then $h(x)\in \Omega$ since for all $j\in\Z^2$, \eqref{eqn_red}
gives
\begin{align*}
\sum_{i\in T} w(i) h(x)|_{T+j}(i) &= \sum_{i \in T}w(i)h(x)(i+j)\\
&= \sum_{i\in T}w(i)x(i+j,i+j)(0)\\
&=\sum_{i\in T}w(i)x(i,i)(j),
\end{align*}
which is $0 \pmod{q}$ since $x\in \Lambda^\Delta$.

To see that $h:\Lambda^\Delta\to \Omega$ is a homeomorphism it suffices
to show $h$ is a continuous bijection. For $f\in \Omega$ define a
function $k(f):\Delta\to\Lambda$ by
\[
k(f)(m,n)=f|_{T(n-m)+m} \text{ for }m\leq n.
\]
Condition \eqref{defX} implies that $k(f)(m,n)$ is a path in
$\Lambda$ of degree $n-m$, so $k$ is degree-preserving. An
application of Proposition~\ref{prop_pathcomposition} says that
$k(f)(m,p)=f|_{T(p-m)+m}$ factors as
\[f|_{T(p-m)+m}=f|_{T(n-m)+m}f|_{T(p-n)+n},\]
which gives
\[k(f)(m,p)=k(f)(m,n)k(f)(n,p),\]
so $k(f)$ is a functor.

We claim that $k:\Omega\to\Lambda^\Delta$ is
the inverse of $h$. We have $h(k(f))=f$ since for $i\in\Z^2$
\[h(k(f))(i) = k(f)(i,i)(0) = f|_{T(0)+i}(0) = f(i+0)=f(i).\]
Now we must show $k(h(x))=x$. Let $i\in T(n-m)$. Then $i=j+l$ for some $j\in T$ and $0\leq l\leq n-m$. We have
\begin{align*}
(k(h(x))(m,n))(i) &= h(x)|_{T(n-m)+m}(i) \\
&= h(x)(i+m) \\
&= h(x)(j+l+m) \\
&= x(j+l+m,j+l+m)(0) \\
&= x(l+m,l+m)(j)\ \text{ by } \eqref{eqn_red}\\
&= x(m,n)|_{T+l}(j) \\
&= x(m,n)(j+l) \\
&= x(m,n)(i),
\end{align*}
so $k(h(x))=x$.

To see that $h$ is continuous, suppose $x_\gamma\to x$ in
$\Lambda^\Delta$. Since $\Omega$ has the product topology, it
suffices to prove that $h(x_\gamma)(i)\to h(x)(i)$ for all $i\in\Z^2$.
Let $i\in\Z^2$. Then $Z(x(i,i),0)$ is an open neighbourhood of $x$
in $\Lambda^\Delta$, so for large $\gamma$, $x_\gamma\in Z(x(i,i),0)$. But
then for large $\gamma$ we have $h(x_\gamma)(i)=x_\gamma(i,i)(0)=x(i,i)(0)=h(x)(i)$, so certainly
$h(x_\gamma(i))\to h(x)(i)$.

For the last part we have $h(\sigma_p(x))=\alpha_p(h(x))$ since
\[h(\sigma_p(x))(i) = (\sigma_px)(i,i)(0) = x(i+p,i+p)(0) = h(x)(i+p) = \alpha_p(h(x))(i).\qedhere\]
\end{proof}

\begin{remark}
Theorem~\ref{thm_infpathsphomeo} implies in particular that the
shift space $\Omega$ associated to the Ledrappier graph
$L(\sock{0.06}{}{}{})$ is the $2$-dimensional Markov system known as
Ledrappier's example (see \cite[Examples 1.8, 2.4]{LS} and
\cite{Ledrapp}).
\end{remark}

\begin{remark}\label{onesided}
There is a one-sided version of Theorem~\ref{thm_infpathsphomeo}. The space
\[
\Omega^+:=\Big\{ f: \N^2\to\Zmod{q} : \sum_{i\in T} w(i) f|_{T+n}(i) = 0 \pmod{q} \text{ for all } n\in\mathbb{N}^2 \Big\} .
\]
has a natural action of $\N^2$, and the $\Z^2$ action on
$\Lambda^\Delta$ restricts to an $\N^2$ action on the one-sided
infinite path space $\Lambda^\infty$. Then the argument of
Theorem~\ref{thm_infpathsphomeo} gives a homeomorphism of
$\Lambda^\infty$ onto $\Omega^+$ which commutes with the actions of
$\N^2$.
\end{remark}

\begin{remark}\label{whyinvert2}
We saw in Remarks~\ref{whyinvcorners} and~\ref{whynot3d} that
relaxing our hypotheses on the rule or using higher-dimensional
tiles would lead to situations where we have to nominate blue-red to
red-blue factorisations to define a $k$-graph $\Lambda$. In the
two-dimensional case, this would mean that if $d(\lambda)=e_1+e_2$,
then $\lambda((c_2+1)e_2)$ will depend on the choice of
$\lambda((c_1+1)e_1)$ as well as the values of $\lambda$ on $T\cup
(T+e_1+e_2)$. So the homeomorphism of Remark~\ref{onesided} will
carry the infinite path space of $\Lambda$ onto a proper subspace of
the shift space $\Omega^+$.
\end{remark}

\section{Aperiodicity}\label{s_aperiodicity}

Aperiodicity is the property of a $k$-graph $\Lambda$ which ensures
that $\Lambda$ has a Cuntz-Krieger uniqueness theorem which says
that all Cuntz-Krieger $\Lambda$-families generate isomorphic
$C^*$-algebras. We will use a formulation of aperiodicity due to
Robertson and Sims: $\Lambda$ is \emph{aperiodic} if for every $v\in
\Lambda^0$ and $m,n\in \N^2$ with $m\neq n$, there is a path
$\lambda\in\Lambda$ satisfying $r(\lambda)=v$, $d(\lambda) \geq
m\vee n$ and
\begin{equation}
\lambda(m,m+d(\lambda)-(m\vee n)) \neq \lambda(n,n+d(\lambda)-(m\vee
n)). \label{eqn_ape}
\end{equation}
It is shown in \cite[Lemma~3.2]{RobS} that this is equivalent to the
aperiodicity hypotheses used in \cite{KP} and \cite{RSY03} (which
phrase aperiodicity as properties of the shifts on the one-sided
path space $\Lambda^\infty$).

To prove aperiodicity of our $2$-graphs we need to make another
restriction on the rule $w$. We say that $w$ \emph{has three
invertible corners} if $w(0)$, $w(c_1e_1)$ and $w(c_2e_2)$ are all
invertible in $\Zmod{q}$ (implicitly demanding that $c_1\geq 1$ and
$c_2\geq 1$). We show in Example~\ref{W(0)not0} that aperiodicity
may fail if $w(0)$ is not invertible. We will also simplify things
by assuming that the trace $t$ is zero, and we will discuss this
hypothesis after the proof of Theorem~\ref{thm_aperiodicity}.

\begin{example}\label{W(0)not0}
Consider the data consisting of the sock tile $T=\{0,e_1,e_2\}$,
$q=2$  $t=0$ and rule defined by $w(0)=0$, $w(e_1)=w(e_2)=1$. The
vertices are:
\[\sock{0.25}{0}{0}{0}\hspace{1cm} \sock{0.25}{1}{0}{0} \hspace{1cm} \sock{0.25}{0}{1}{1} \hspace{1cm} \sock{0.25}{1}{1}{1}\;.\]
Since every vertex $v$ has $v(e_1)=v(e_2)$, every path is constant
along the short diagonals $n_1+n_2=c$. In other words, for every
path $\lambda$ and every $n$, we have
$\lambda(n)=\lambda(n+e_1-e_2)$ whenever $n$ and $n+e_1-e_2$ lie in
the domain of $\lambda$. This implies in particular that for every
path $\lambda$ with $d(\lambda)\geq (1,1)$, we have
\[
\lambda(e_2, d(\lambda)-e_1)=\lambda(e_1,d(\lambda)-e_2),
\]
so \eqref{eqn_ape} fails for every $v$ with $m=e_2$ and $n=e_1$.
\end{example}

\begin{theorem}\label{thm_aperiodicity}
If the rule $w$ in the basic data $(T,q,0,w)$ has $c_1\geq 1$,
$c_2\geq 1$ and three invertible corners, then the associated
$2$-graph $\Lambda$ is aperiodic.
\end{theorem}

For the proof of the theorem we need to know that $\Lambda$ is
\emph{strongly connected} in the sense that every $v\Lambda^*u$ is
non-empty.

\begin{prop}\label{prop_transitivity}
Suppose that $k\in\N$ satisfies $(k-1)(e_1+e_2)\in T$ and
$k(e_1+e_2)\notin T$. Then for every $v,u\in\Lambda^0$ there exists $\lambda\in\Lambda^{k(e_1+e_2)}$ such that $r(\lambda)=v$ and
$s(\lambda)=u$.
\end{prop}

The proof of Proposition~\ref{prop_transitivity} depends on the
following variant of Proposition~\ref{prop_skeleton}(b).

\begin{lemma}\label{lem_existpath11}
If $v,u\in\Lambda^0$ satisfy
\begin{equation}\label{eqn_oneoneoverlap}
v(j) = u(j-e_1-e_2) \text{ for } j\in T\cap (T+e_1+e_2),
\end{equation}
then there is a unique path $\mu\in\Lambda^{e_1+e_2}$ such that $r(\mu)=v$ and
$s(\mu)=u$.
\end{lemma}

\begin{proof}We define
\[
\mu(j)=\begin{cases}
v(j)&\text{if $j\in T$}\\
u(j-(e_1+e_2))&\text{if $j\in (T+e_1+e_2)\backslash T$.}
\end{cases}
\]
then we obviously have $\mu|_T=v$, and \eqref{eqn_oneoneoverlap}
says that $\mu|_{T+e_1+e_2}=u$. Now we observe that because the
corners $w(c_ie_i)$ are invertible, there are unique values of
$\mu((c_i+1)e_i)$ such that $\mu|_{T+e_i}$ belongs to $\Lambda^0$.
So there is exactly one path $\mu$ with the required property.
\end{proof}

\begin{proof}[Proof of Proposition~\ref{prop_transitivity}]
We will prove by induction on $p$ that for $0\leq p\leq k$, there
exists $\mu^p\in \Lambda^{p(e_1+e_2)}$ such that $r(\mu^p)=v$ and
\begin{equation}
s(\mu^p)(j)=u(j-(k-p)(e_1+e_2)) \text{ for }
j\in T\cap(T+(k-p)(e_1+e_2)). \label{eqn_ztarget}
\end{equation}
Then $\mu:=\mu^k$ is the required path.

For $p=0$, we take $\mu^0:=v$. Suppose that $0\leq p<k$ and we have
$\mu^p$ with the required properties. Now we define
\begin{equation}\label{defvp+1}
v^{p+1}(i)=\begin{cases}
s(\mu^p)(i+e_1+e_2)&\text{for $i\in T\cap(T-e_1-e_2)$}\\
u(i-(k-p-1)(e_1+e_2))&\text{for $i\in T\cap(T+(k-p-1)(e_1+e_2))$;}
\end{cases}
\end{equation}
if $j$ belongs to both sets on the right-hand side, then we can
apply \eqref{eqn_ztarget} with $j=i+e_1+e_2$ and deduce that the two
possible values for $v^{p+1}(i)$ coincide. We now define
$v^{p+1}(i)$ arbitrarily for other points $i$ in $T\backslash
\{c_1e_1\}$, and set
\[
v^{p+1}(c_1e_1):=w(c_1e_1)^{-1}\Big(\sum_{i\in T\backslash \{c_1e_1\}}w(i)v^{p+1}(i)\Big),
\]
so that $v^{p+1}\in \Lambda^0$. The first option in \eqref{defvp+1}
implies that the pair $s(\mu^p)$ and $v^{p+1}$ satisfy
\eqref{eqn_oneoneoverlap}, and hence by Lemma~\ref{lem_existpath11}
there exists a path $\nu\in \Lambda^{e_1+e_2}$ with
$r(\nu)=s(\mu^p)$ and $s(\nu)=v^{p+1}$. Now we take $\mu^{p+1}$ to
be the composition $\mu^p\nu$, and the second option in
\eqref{defvp+1} implies that $s(\mu^{p+1})$ satisfies
\eqref{eqn_ztarget}.
\end{proof}

\begin{proof}[Proof of Theorem~\ref{thm_aperiodicity}]
We fix $v\in \Lambda^0$ and $m,n\in \N^2$ with $m\not= n$. We choose
a path $\mu$ with $r(\mu)=v$ and $d(\mu)=m\vee n$. We aim to extend
$\mu$ to a path $\lambda$ satisfying \eqref{eqn_ape}. If the
vertices $\mu|_{T+m}$ and $\mu|_{T+n}$ are different, then
$\lambda:=\mu$ will do. So we suppose that $\mu|_{T+m}=\mu|_{T+n}$.
We deal separately with the cases where $m$ and $n$ are comparable
in the sense that either $m\leq n$ or $n\leq m$, and where they are
not comparable.

Suppose first that $m$ and $n$ are comparable, say $m\leq n$. Since
$m\not=n$, there exists $i$ such that $m+e_i\leq n$, and then we
have
\[
r(\mu(m,m+e_i))=\mu|_{T+m}=\mu|_{T+n}=s(\mu).
\]
Since $c_i\geq 1$, Proposition~\ref{prop_skeleton}(c) implies that
$|s(\mu)\Lambda^{e_i}|>1$ and so there exists $\nu\in
s(\mu)\Lambda^{e_i}$ such that $\nu\neq \mu(m,m+e_i)$. Then
$\lambda:=\mu\nu$ has the required properties:
\[
\lambda(m,m+d(\lambda)-(m\vee n))= \lambda(m,m+e_i)= \mu(m,m+e_i)
\]
is not equal to
\[
\lambda(n,n+d(\lambda)-(m\vee n))= \lambda(n,n+e_i)= \nu.
\]

Now suppose that $m$ and $n$ are not comparable, say $m_1>n_1$ and
$m_2<n_2$. This is where we use the extra hypotheses on the rule $w$
and the trace $t$. Since $t=0$, the identically zero function
$v_0:T\to \Z/q\Z$ defines a vertex $v_0$, and the identically zero
function $x:\N^2\to \Z/q\Z$ defines an infinite path
$x\in\Lambda^{\infty}$ (via the homeomorphism of
Remark~\ref{onesided}). Since $\Lambda^{e_1}v_0$ has more than one
element, and there is just one blue edge from $v_0$ to $v_0$ (see
Proposition~\ref{prop_skeleton}), there must be a blue edge $\beta$
with $s(\beta)=v_0$ and $r(\beta)\not= v_0$. By
Proposition~\ref{prop_transitivity}, there is a path $\alpha$ with
$r(\alpha)=s(\mu)$ and $s(\alpha)=r(\beta)$. We claim that
\[
\lambda:=\mu\alpha\beta x(0,(m\vee n)-(m\wedge n)-e_1)=\mu\alpha\nu,\ \text{ say,}
\]
satisfies \eqref{eqn_ape}; indeed, we claim that the two paths in
\eqref{eqn_ape} have different sources. Since
$d(\lambda)=d(\mu\alpha)+(m\vee n) -(m\wedge n)$ and
\begin{align*}
\lambda|_{T+m+d(\lambda)-(m\vee n)}&=
\lambda|_{T+m+d(\mu\alpha)-(m\wedge n)}\\
&= \nu|_{T+m-(m\wedge n)}\\
&=x|_{T+(m_1-n_1-1)e_1}\\
&=v_0,
\end{align*}
it suffices to show that
\begin{align*}
\lambda|_{T+n+d(\lambda)-(m\vee n)}
= \nu|_{T+n-(m\wedge n)}= \nu|_{T+(n_2-m_2)e_2}
\end{align*}
is not equal to $v_0$.

We suppose that there exists $p\in \N$ such that
$\nu|_{T+pe_2}=v_0$, and look for a contradiction. Then there is a
smallest such $p$, and since $\nu|_T=r(\beta)\not=v_0$, we then have
$p>0$. Now $\nu|_{T+(p-1)e_2}\in \Lambda^0$ implies
\begin{equation}\label{nuincol1}
w(0)\nu((p-1)e_2)=-\sum_{i\in T\backslash\{0\}} w(i)\nu(i+(p-1)e_2);
\end{equation}
since we have $\nu(l)=x(l)=0$ whenever $l_1>0$, \eqref{nuincol1}
implies that
\begin{equation}
w(0)\nu((p-1)e_2)=-\sum_{j=1}^{c_2} w(je_2)\nu((j+p-1)e_2),
\end{equation}
which is $0$ because $\nu((k+p)e_2)=\nu|_{T+pe_2}(ke_2)=v_0(ke_2)=0$
for $k\geq 0$. Since $w(0)$ is invertible, this implies that
$\nu((p-1)e_2)=0$. Thus we have $\nu|_{T+(p-1)e_2}=v_0$, and this
contradicts the choice of $p$. Thus for every $p$, $\nu|_{T+pe_2}$
is not equal to $v_0$, and in particular $\nu|_{T+(n_2-m_2)e_2}$ is
not equal to $v_0$, as required.
\end{proof}

\begin{remark}\label{remark_isomorphictotrace0}
The preceding proof also works when $t\not=0$ provided there is a
vertex $v_0$ which is constant, say $v_0(m)=c$ for all $m\in T$.
There is such a vertex if and only if there exists $c\in \Z/q\Z$
such that
\begin{equation}
c\Big(\sum_{i\in T}w(i)\Big) = t \pmod{q}. \label{eqn_loopbothcolourscondition}
\end{equation}
However, we do not obtain any new $2$-graphs this way: if there is
such a $c$, then $\Lambda(T,q,t,w)$ is isomorphic to
$\Lambda(T,q,0,w)$. To see this, note that for every path $\lambda$
in $\Lambda(T,q,0,w)$, $\lambda_t:i\mapsto \lambda(i)+c\pmod{q}$ is
a path in $\Lambda(T,q,t,w)$, and the map $\lambda\mapsto \lambda_t$
is an isomorphism of $\Lambda(T,q,0,w)$ onto $\Lambda(T,q,t,w)$.

It is easy to find examples where
\eqref{eqn_loopbothcolourscondition} has no solution $c$. For
example, if $|T|$ is even, $q=2$, $t=1$ and $w\equiv 1$, we have
$\sum w(i)=|T|$ and  $c|T|=1\pmod{2}$ has no solutions. We do not
have general criteria for aperiodicity when
\eqref{eqn_loopbothcolourscondition} has no solution.
\end{remark}

\section{The $C^*$-algebras}\label{s_cstaralgs}

We now summarise the properties of the $C^*$-algebras of the
$2$-graphs $\Lambda(T,q,0,w)$.

\begin{theorem}\label{thm_spin}
Suppose $(T,q,0,w)$ is basic data with $c_1\geq 1$ and $c_2\geq 1$,
and the rule $w$ has three invertible corners.  Then
$C^*(\Lambda(T,q,0,w))$ is unital, nuclear, simple and purely
infinite, and belongs to the bootstrap class~$\mathcal{N}$.
\end{theorem}

\begin{proof}
We write $\Lambda$ for $\Lambda(T,q,0,w)$. We begin by observing
that $C^*(\Lambda)$ is unital because $\Lambda^0$ is finite, and is
nuclear and belongs to the bootstrap class by
\cite[Theorem~5.5]{KP}. It follows easily from
Proposition~\ref{prop_transitivity} that $\Lambda$ is cofinal: if
$x\in \Lambda^\infty$ and $v\in \Lambda^0$, then there is a path
from $r(x)$ to $v$.  Since we know from
Theorem~\ref{thm_aperiodicity} that $\Lambda$ is aperiodic (that is,
satisfies property (iv) of \cite[Lemma~3.2]{RobS}), it follows from
Theorem~3.1 and Lemma~3.2 of \cite{RobS} that $C^*(\Lambda)$ is
simple.

To see that $\Lambda$ is purely infinite, we need to check that
every vertex $v$ can be reached from a ``loop with an entrance''
(see \cite[Proposition~8.8]{Sims}). But we know that $v$ receives at
least two blue edges $\alpha$, $\beta$, and then
Proposition~\ref{prop_transitivity} implies that there is a path
$\nu$ from $v$ to $s(\alpha)$, so there is a path $\mu=\alpha\nu$
with $d(\mu)\not=0$ such that $r(\mu)=s(\mu)=v$. Since $\beta$ is an
entrance to $\mu$, we have verified the hypothesis of
\cite[Proposition~8.8]{Sims}, and can deduce that $C^*(\Lambda)$ is
purely infinite.
\end{proof}

\begin{remark}
We have appealed to \cite[Proposition~8.8]{Sims} rather than
\cite[Proposition~4.9]{KP} because the latter is not correct as it
stands. For example, the $2$-graphs in \cite[Figures~3~and~4]{PRRS}
satisfy the hypothesis of \cite[Proposition~4.9]{KP}, but their
$C^*$-algebras are A$\T$-algebras and hence are not purely infinite.
\end{remark}

\section{K-theory}\label{s_kthry}

Theorem~\ref{thm_spin} implies that, when the rule has three
invertible corners, the $C^*$-algebra falls into the class which is
classified by the celebrated theorem of Kirchberg and Phillips,
which says that $C^*(\Lambda)$ is determined up to isomorphism by
its $K$-theory \cite{K2, P, Ror}. So we want to compute the
$K$-groups of $C^*(\Lambda)$.

Suppose we have basic data satisfying the hypotheses of
Proposition~\ref{prop_skeleton}, so that in particular the associated $2$-graph
$\Lambda$ is finite with no sources, and the methods of \cite{E} apply.  Let $B$
and $R$ be the vertex matrices of $\Lambda$, defined for
$u,v\in\Lambda^0$ by
\begin{align*}
B(u,v)&=\#\{\lambda\in \Lambda^{e_1}:r(\lambda)=u,s(\lambda)=v\}\\
R(u,v)&=\#\{\lambda\in \Lambda^{e_2}:r(\lambda)=u,s(\lambda)=v\};
\end{align*}
the matrices $B$ and $R$ are the vertex matrices of the \emph{blue
graph} $B\Lambda:=(\Lambda^0,\Lambda^{e_1},r,s)$ and the \emph{red
graph} $R\Lambda:=(\Lambda^0,\Lambda^{e_2},r,s)$. The entries
$BR(u,v)$ in the product $BR$ are the numbers of blue-red paths from
$v$ to $u$, which the factorisation property implies are the same as
the entries $RB(u,v)$ in $RB$; in other words, $BR=RB$. Let
$\delta_1:\ZL0\oplus \ZL0\to \ZL0$ and $\delta_2:\ZL0\to \ZL0\oplus
\ZL0$ be the maps with matrices
\[
\delta_1=\begin{pmatrix}1-B^t &
1-R^t \end{pmatrix} \text{ and } \delta_2=\begin{pmatrix}R^t-1 \\
1-B^t \end{pmatrix}.
\]
Then Proposition~3.16 of \cite{E} says that the
$K$-groups are given by
\begin{align*}
K_0(\cstarlambda) &\cong \coker \delta_1\oplus \ker \delta_2 \\
K_1(\cstarlambda) &\cong \ker \delta_1/\img \delta_2.
\end{align*}

We were able to calculate the size of the $K$-groups for a large
number of examples by implementing the following procedure in the
\texttt{Magma} computational algebra system. \texttt{Magma} recognises that we are
dealing with integer matrices and so it performs calculations over
the integers; for example, on being asked to find a basis for the
columnspace of an integer matrix it returns an integer basis. When
calculating $|K_0(\cstar(\Lambda))|$, we noticed that
$\ker\delta_2=0$ in every example. To calculate $|\coker\delta_1|$,
we find a basis matrix $M$ whose columns are an integer basis for
the columnspace of the matrix of $\delta_1$. Then
\[|K_0(\cstarlambda)|=|\coker\delta_1|=|\det M|.\]
To calculate $|K_1(\cstar(\Lambda))|$, first we find a basis matrix
$H$ for $\ker\delta_1$. Since the columns of $H$ are linearly
independent, for each column vector $z$ of the matrix of $\delta_2$
the equation $Hw=z$ has a unique solution $w$. Form the matrix whose
columns are the solutions $w$; then the $i$th column of $W$ contains
the co-ordinates of the basis vector $z$ with respect to the basis
for $\ker\delta_1$. Thus
\[|K_1(\cstarlambda)|=|\ker\delta_1/ \img\delta_2|=|\det W|.\]
We give details of these calculations for the Ledrappier graph.

\begin{example}\label{ex_sockKthry}
Consider the sock tile $\sock{0.1}{}{}{}$ with $q=2$, $t=0$ and rule
$w\equiv 1$. Vertex matrices for the $2$-graph $\Lambda$ are
\[B=\begin{pmatrix}
  1 & 1 & 0 & 0 \\
  0 & 0 & 1 & 1 \\
  1 & 1 & 0 & 0 \\
  0 & 0 & 1 & 1 \\
\end{pmatrix}\text{ and }
R=\begin{pmatrix}
  1 & 1 & 0 & 0 \\
  0 & 0 & 1 & 1 \\
  0 & 0 & 1 & 1 \\
  1 & 1 & 0 & 0 \\
\end{pmatrix}.\]
The matrices $\delta_1:\Z^8\to\Z^4$ and
$\delta_2:\Z^4\to\Z^8$ are
\[\delta_1=
\begin{pmatrix}
  0 & 0 & -1 & 0 & 0 & 0 & 0 & -1 \\
  -1 & 1 & -1 & 0 & -1 & 1& 0  & 1 \\
  0 & -1 & 1& -1 & 0 & -1 & 0 & 0 \\
  0 & -1& 0 & 0 & 0 & -1 & -1 & 1 \\
\end{pmatrix}
\text{ and } \delta_2=
\begin{pmatrix}
  0 & 0 & 0 & 1 \\
  1 & -1 & 0 & 1 \\
  0 & 1 & 0 & 0 \\
  0 & 1 & 1 & -1 \\
  0 & 0 & -1 & 0 \\
  -1 & 1 & -1 & 0 \\
  0 & -1 & 1 & -1 \\
  0 & -1 & 0 & 0 \\
\end{pmatrix}.\]
Here, $\delta_1$ is onto so we can choose a basis for $\img\delta_1$
such that $M$ is the $4\times 4$ identity matrix; hence $|\det M|=1$, and $K_0(\cstarlambda)=0$. \texttt{Magma} gives us the matrices $H$
and $W$ below, which satisfy $HW=\delta_2$.
\[H=\begin{pmatrix}
      1 & 0 & 0 & 0 \\
      0 & 1 & 0 & 0 \\
      0 & 0 & 1 & 0 \\
      0 & 0 & 0 & 1 \\
      -1 & 0 & 1 & -1 \\
      0 & -1 & 1 & -1 \\
      0 & 0 & -2 & 1 \\
      0 & 0 & -1 & 0 \\
    \end{pmatrix}
\text{ and }W=
    \begin{pmatrix}
      0 & 0 & 0 & 1 \\
      1 & -1 & 0 & 1 \\
      0 & 1 & 0 & 0 \\
      0 & 1 & 1 & -1 \\
    \end{pmatrix}.
\]
Then since $|\det W|=1$, we have $K_1(\cstarlambda)=0$.
\end{example}

Some of the results of our calculations are listed in
Table~\ref{table_kthrytable}. We begin with some explanatory comments.
\begin{itemize}
\item A tile $T$ can be uniquely described by the lengths of its rows from longest to shortest.
For example, in the table we write $[2,1]$ for the sock.
\item The tile obtained by reflecting $T$ about the line $y=x$ is called the \emph{conjugate tile} of $T$.
For example, the conjugate of the tile $[3,1]$ is $[2,1,1]$ and the
sock tile is its own conjugate. A tile and its conjugate give $\cstar$-algebras with the
same $K$-theory since this amounts to swapping the roles of $B$ and
$R$ in the $K$-theory formulas. So in the table we list only one out of
each pair of conjugate tiles.
\item The results in the table refer to basic data with $t=0$ and $w\equiv 1$.
We also performed calculations for other rules and traces, but we obtained the same values of $|K_0|$ and $|K_1|$. A partial explanation for
this is in Remark~\ref{remark_isomorphictotrace0}.
\item Blank spaces in the table would require calculations beyond sensible computation time. We were able to do more calculations when $q=2$, but the results did not reveal any interesting new phenomena.
\end{itemize}

\begin{table}
\begin{center}\begin{tabular}{|l|l|cc|cc|cc|cc|}
  \hline
  size & tile & $q=2$ &   & $q=3$ & & $q=4$ &  & $q=5$ &  \\
  && $|K_0|$ & $|K_1|$ &  $|K_0|$ & $|K_1|$ & $|K_0|$ & $|K_1|$ & $|K_0|$ & $|K_1|$  \\
  \hline
  3 & [3] & 3 & 3 & 8 & 8 & 15 & 15 & 24 & 24 \\
    & [2,1] & 1 & 1 & 2 & 2 & 3 & 3 & 4 & 4 \\ \hline
  4 & [4] & 7 & 7 & 26 & 26 & 63 & 63 & 124 & 124 \\
    & [3,1] & 1 & 1 & 2 & 2 & 3 & 3 & 4 & 4 \\
    & [2,2] & 1 & 1 & 2 & 2 & 3 & 3 & 4 & 4 \\ \hline
  5 & [5] & 15 & 15 & 80 & 80 & 255 & 255 & 624 & 624 \\
    & [4,1] & 1 & 1 & 2 & 2 & 3 & 3 & 4 & 4 \\
    & [3,2] & 1 & 1 & 2 & 2 & 3 & 3 & 4 & 4 \\
    & [3,1,1] & 3 & 3 & 8 & 8 & 15 & 15 & 24 & 24 \\ \hline
  6 & [6] & 31 & 31 & 242 & 242 & 1023 & 1023 &  &  \\
    & [5,1] & 1 & 1 & 2 & 2 & 3 & 3 &  &  \\
    & [4,2] & 1 & 1 & 2 & 2 & 3 & 3 &  &  \\
    & [4,1,1] & 1 & 1 & 2 & 2 & 3 & 3 &  &  \\
    & [3,3] & 1 & 1 & 2 & 2 & 3 & 3 &  &  \\
    & [3,2,1] & 3 & 3 & 8 & 8 & 15 & 15 &  &  \\ \hline
  7 & [7] & 63 & 63 & 728 & 728 &  &  &  &  \\
    & [6,1] & 1 & 1 & 2 & 2 &  &  &  &  \\
    & [5,2] & 1 & 1 & 2 & 2 &  &  &  &  \\
    & [5,1,1] & 3 & 3 & 8 & 8 &  &  &  &  \\
    & [4,3] & 1 & 1 & 2 & 2 &  &  &  &  \\
    & [4,2,1] & 1 & 1 & 2 & 2 &  &  &  &  \\
    & [4,1,1,1] & 7 & 7 & 26 & 26 &  &  &  &  \\
    & [3,3,1] & 3 & 3 & 8 & 8 &  &  &  &  \\ \hline
  \end{tabular}\end{center}
\caption{Table of $K$-theory calculations} \label{table_kthrytable}
\end{table}

A more detailed look at the results of our calculations suggests the following

\begin{conjecture}
\begin{enumerate}
\item $\ker\delta_2=0$.
\smallskip
\item $|K_0(\cstar(\Lambda))|=|K_1(\cstar(\Lambda))|$.
\smallskip
\item $|K_i(\cstar(\Lambda)))|$ always has the form $q^n-1$. Our calculations are consistent with the formula $|K_i(\cstar(\Lambda))|=(q^{c_2}-1, q^{c_1}-1)$.
\end{enumerate}
\end{conjecture}
We prove Conjectures $(1)$ and $(2)$ in Theorems
\ref{thm_kerdelta2trivial} and \ref{thm_K0eqK1} below. We do not
know whether $K_0(\cstar(\Lambda))$ is isomorphic to
$K_1(\cstar(\Lambda))$ in general, though calculations in \texttt{Magma} confirm that $K_i(\cstar(\Lambda))$ is cyclic in all the examples listed in Table~\ref{table_kthrytable}, and hence $K_0$ is isomorphic to $K_1$ for all these examples. However, this is automatically true in most cases because there is only one group of the given order, so the number of examples we have considered where there is something to prove (that is, the ones where $|K_i|=4$ or $8$) is fairly small. Certainly we have not yet identified a potential reason for the existence of such an isomorphism, and our proof of Conjecture (2) does not help. We have only numerical evidence for Conjecture $(3)$.

\subsection*{Implications for the classification}

The graphs whose $K$-theory is computed in
Table~\ref{table_kthrytable} are of two types. For tiles with
$c_2=0$, the blue graph consists of disjoint cycles, and the
aperiodicity condition of \cite{RobS} fails, so the
$\cstar$-algebras of these graphs are not simple. (Their structure
is nevertheless quite intricate and will be discussed in a future
paper.) For all other graphs, the basic data satisfies the
hypotheses of Theorem~\ref{thm_spin}, and hence the
$\cstar$-algebras are simple and satisfy the hypotheses of the
Kirchberg-Phillips Theorem. The Kirchberg-Phillips Theorem (as stated
in \cite[Theorem~8.4.1(iv)]{Ror}, for example) says that two suitable unital
$\cstar$-algebras $A$ and $B$ are isomorphic if and only if
$K_1(A)\cong K_1(B)$ and there is an isomorphism of $K_0(A)$ onto
$K_0(B)$ which takes the class $[1_A]$ of the identity to
$[1_B]$. When $|K_0(\cstarlambda)|=|K_1(\cstarlambda)|=1$, the last
condition is trivially satisfied and $\cstarlambda$ is isomorphic to
the Cuntz algebra $\O_2$. (Somewhat disappointingly, the Ledrappier
graph is one of these graphs.) When $|K_i|>1$, we computed the class
of $[1]=\sum_{v\in\Lambda^0}[p_v]$ in
$K_0(\cstar(\Lambda))=\coker\delta_1$ (see \cite[Corollory~5.1]{E}),
and found that it is always a generator for $K_0(\cstarlambda)$. So
in all our examples, $K_0(\cstarlambda)$ is cyclic. We do not
know whether this is always true. To sum up: if $\Lambda_1$ and
$\Lambda_2$ are any graphs in the table with $c_2\geq 1$, and if
$K_0(\cstar(\Lambda_1))=K_0(\cstar(\Lambda_2))$, then
$\cstar(\Lambda_1)\cong \cstar(\Lambda_2)$.

None of the $\cstar$-algebras of graphs in
Table~\ref{table_kthrytable} with non-zero $K$-theory can be
isomorphic to the $\cstar$-algebra of an ordinary directed graph
$E$, because $K_1(\cstar(E))$ is always free (being a subgroup of
the free group $\Z^{E^0}$).

\section{$K$-theory results}\label{s_kthryresults}

Let $T$ be a tile, and again write $(c_1,c_2)=\bigvee\{j:j\in T\}$.
For $0\leq i \leq c_1$, we let $h_i$ denote the second coordinate of
the top box $(i,h_i)$ in each column of $T$; for $0\leq
i\leq c_2$, $w_i$ is the first coordinate of the
right-hand box $(w_i,i)$ in each row.

In this section we prove conjectures $(1)$ and $(2)$ about
$K_*(\cstar(\Lambda(T,q,t,w)))$, under some mild hypotheses on the shape of the tile.

\begin{theorem}\label{thm_kerdelta2trivial}
Suppose we have basic data $(T,q,t,w)$ in which $w$ has invertible
corners and $c_1,c_2\geq 1$. Suppose further that  either
$h_0>h_1$ or $w_0>w_1$. If $B$ and $R$ are the vertex matrices
associated to $\Lambda=\Lambda(T,q,t,w)$, then the map
\[
\delta_2=\begin{pmatrix}R^t-1\\1-B^t\end{pmatrix}:\ZL0\to\ZL0\oplus\ZL0
\]
has trivial kernel and $K_0(\cstarlambda)=\coker \delta_1$.
\end{theorem}

Proposition~3.16 of \cite{E} says that the bijection of
$\coker\delta_1$ onto $K_0(\cstarlambda)$ carries the generator
$\delta_v$ of $\ZL0$ into $[p_v]$, and therefore
Theorem~\ref{thm_kerdelta2trivial} says that these generate
$K_0(\cstarlambda)$. The image of $\delta_1$ is then generated by
the images of the elements $(1-B^t)\delta_v$ and $(1-R^t)\delta_v$.
Thus Theorem~\ref{thm_kerdelta2trivial} says that, for one of our
$2$-graphs $\Lambda$, $K_0(\cstarlambda)$ is generated by
$\{[p_v]:v\in\Lambda^0\}$ modulo the relations
\[[p_v]=\sum_{r(e)=v,d(e)=e_1}[p_{s(e)}],\hspace{1cm} [p_v]=\sum_{r(e)=v,d(e)=e_2}[p_{s(e)}]\]
imposed by the blue and red Cuntz-Krieger relations.

Theorem~\ref{thm_kerdelta2trivial} will follow immediately from the following proposition. For the rest of the section, we fix a set of basic data $(T,q,t,w)$ in which $w$ has invertible corners and $c_1,c_2\geq 1$.

\begin{prop}\label{prop_kernel0}
Let $B$ and $R$ be the
vertex matrices of $\Lambda=\Lambda(T,q,t,w)$.
\begin{enumerate}
\item If $h_0>h_1$ then the map $1-B^t:\ZL0\to\ZL0$ has trivial
kernel.
\item If $w_0>w_1$ then the map $1-R^t:\ZL0\to\ZL0$ has trivial kernel.
\end{enumerate}
\end{prop}

To prove this we use the special structure of the vertex matrices of
$\Lambda$. By symmetry it suffices to prove part $(1)$. From
Proposition~\ref{prop_skeleton} we know that $B$ is a
$\{0,1\}$-matrix; that the number of $1$s in each row/column is
$q^{c_2}$; and that any two rows/columns are either equal or
orthogonal. The crucial observation is that the matrices with these
properties are the ones which arise as the vertex matrices of dual
graphs. Recall from \cite[page~17]{R}, for example, that the
\emph{dual graph} of a directed graph $E$ is the directed graph
$\widehat{E}$ with $\widehat E^0:=E^1$, $\widehat E^1:=\{(e,f)\in
E^1\times E^1:r(f)=s(e)\}$, and range and source maps given by
$r(e,f)=r(e)$ and $s(e,f)=s(f)$.

To describe the graphs whose duals arise we need some notation. Let
$S$ be the tile $S:=T\cap(T-e_1)$ and let $S^+$ be the tile
$S^+:=S\cup\{(h_1+1)e_2\}$. In the visual model, $S$ is the tile
obtained from by deleting the first column and shifting one unit to
the left, and $S^+$ is obtained from $S$ by adding one box to the
top of its first column. For a directed graph $F$ and an integer
$n\geq 1$, the directed graph $nF$ has vertex set $(nF)^0=F^0$, edge
set
\[(nF)^1=F^1\times \{1,\ldots, n\}=\{(f,i): f\in F^1, 1\leq i\leq
n\}\] and range and source maps given by $r(f,i)=r(f)$ and
$s(f,i)=s(f)$. Then the vertex matrix of $nF$ is $n$ times the
vertex matrix of $F$. (Note if $n=1$ then $1F\cong F$.)

\begin{prop}\label{prop_dual}
Suppose that $h_0>h_1$. Set
$r_B=q^{h_0-h_1-1}$, and let $B\Lambda(S^+,q,0,1)$ be the blue
graph of the tile $S^+$ with alphabet $q$, trace $0$ and rule which
is identically $1$. Then the blue graph $B\Lambda$ of
$\Lambda(T,q,t,w)$ is isomorphic to the dual of
$r_BB\Lambda(S^+,q,0,1)$.
\end{prop}

To prove this we need the following lemma.

\begin{lemma}\label{lemma_vtxset}
Let $v_1,v_2\in\Lambda^0$. Then the set
\begin{equation}\label{eqn_setui}
\{u\in \Lambda^0: u|_{S+e_1}=v_2|_S \text{ and } u|_S=v_1|_S\}
\end{equation}
contains $r_B=q^{h_0-h_1-1}$ vertices if $v_1,v_2$ satisfy
\begin{equation}\label{eqn_Soverlap}
v_1(i)=v_2(i-e_1) \text{ for } i\in S\cap (S+e_1),
\end{equation}
and is empty otherwise.
\end{lemma}

\begin{proof}
If \eqref{eqn_Soverlap}, then define a function $u:T\to\Zmod{q}$ by
$u|_{S+e_1}=v_2|_S$ and $u|_S=v_1|_S$. Since
$|T\setminus(S\cup(S+e_1))|=h_0-h_1>0$, there are $r_B$ such
functions $u$ which define vertices in $\Lambda^0$.
\end{proof}

Define a relation $\sim$ on $\Lambda^0$ by
\[v_1\sim v_2 \Longleftrightarrow v_1|_S=v_2|_S.\]
It is straightforward to check that $\sim$ is an equivalence
relation. Let $[v]$ denote the equivalence class of $v\in\Lambda^0$
under $\sim$. By definition of $\sim$ the set in
Lemma~\ref{lemma_vtxset} does not change if we replace $v_1$ and
$v_2$ by other elements of $[v_1]$ and $[v_2]$. So for $v_1,v_2$ satisfying \eqref{eqn_Soverlap}, we can list the
vertices in the set \eqref{eqn_setui} as $u_i([v_1],[v_2])$ for $1\leq i\leq r_B$.

\begin{proof}[Proof of Proposition~\ref{prop_dual}]
Let $F$ be the directed graph with vertices $F^0=\Lambda^0/\sim$ and
edges
\[F^1=\{([v_1],[v_2])\in F^0\times F^0: v_1,v_2 \text{ satisfy } \eqref{eqn_Soverlap}\}.\]
We prove first that $B\Lambda$ is isomorphic to the dual of the
directed graph $r_BF$, and then that $F$ is isomorphic to
$B\Lambda(S^+,q,0,1)$.

By definition, $r_BF$ has vertices $F^0$ and edges
\[(r_BF)^1=\{([v_1],[v_2],i):([v_1],[v_2])\in F^1, 1\leq i\leq r_B\}.\]
So the dual $\widehat{r_BF}$ has vertices
$(\widehat{r_BF})^0=(r_BF)^1$ and there is an edge from
$([v_3],[v_4],j)$ to $([v_1],[v_2],i)$ if and only if $[v_2]=[v_3]$.

Define a map $\phi^0:(\widehat{r_BF})^0\to\Lambda^0$ by
$\phi^0([v_1],[v_2],i)=u_i([v_1],[v_2])$. Then $\phi^0$ is a
bijection since $u_i([v_1],[v_2])$ is uniquely determined by
$([v_1],[v_2])$ and $i$, and every $v\in\Lambda^0$ belongs to a set
in \eqref{eqn_setui} for some $v_1$ and $v_2$ (for example, take
$v_1=v$ and $v_2$ to be any vertex adjacent to $v$).

Suppose there is an edge in $\widehat{r_BF}$ from $([v_3],[v_4],j)$
to $([v_1],[v_2],i)$ --- that is, suppose $[v_2]=[v_3]$. Then
Proposition~\ref{prop_skeleton} says there is a unique edge in
$B\Lambda$ from $\phi([v_3],[v_4],j)=u_j([v_3],[v_4])$ to
$\phi([v_1],[v_2],i)=u_i([v_1],[v_2])$ since
\[ u_i([v_1],[v_2])|_{S+e_1}=v_2|_S=v_3|_S=u_j([v_3],[v_4])|_S. \]
Define the map $\phi^1:(\widehat{r_BF})^1\to\Lambda^{e_1}$ by taking $\phi^1(([v_1],[v_2],i),([v_3],[v_4],j))$ to be the unique edge in
$B\Lambda$ with source $u_j([v_3],[v_4]])$ and range
$u_i([v_1],[v_2])$. Then $\phi^1$ is a bijection since $\phi^0$ is.
We have $r\circ\phi^1=\phi^0\circ r$ since
\begin{align*}
r(\phi(([v_1],[v_2],i),([v_3],[v_4],j))) &= u_i([v_1],[v_2]) \\
&= \phi^0([v_1],[v_2],i) \\
&= \phi^0(r(([v_1],[v_2],i),([v_3],[v_4],j)))
\end{align*}
and similarly $s\circ\phi^1=\phi^0\circ s$. Thus $\phi=(\phi^0,\phi^1)$ is a graph isomorphism from
$\widehat{r_BF}$ to $B\Lambda$.

It remains to show that $F$ is isomorphic to $B\Lambda(S^+,q,0,1)$.
Let $[v]\in F^0$. Define $v^+:S^+\to\Zmod{q}$ by
\[ v^+|_S = v|_S \text{ and } v^+(0,h_1+1) = -\sum_{j\in S}v(j)\pmod{q}.\]
This is well-defined since $S^+\setminus S=\{(0,h_1+1)\}$ and each
element in $[v]$ takes the same values on $S$. We also have that
$v^+$ is uniquely determined by $[v]$, and $v^+$ is clearly a vertex
in $B\Lambda(S^+,q,0,1)$. Define $\psi^0:F^0\to B\Lambda(S^+,q,0,1)^0$ by $\psi^0([v])=v^+$. Then we claim that
$\psi^0$ is a bijection. It is one-to-one because $v^+$ is uniquely determined by $[v]$. To see that $\psi^0$ is onto, let $u$ be a vertex in
$B\Lambda(S^+,q,0,1)$ and suppose $u^-\in\Lambda^0$ with
$u^-|_{S^+}=u$. Then $(u^-)^+|_S=u^-|_S=u|_S$ which implies
$(u^-)^+|_{S^+}=u|_{S^+}$; this says that $\psi([u^-])=(u^-)^+=u$,
so $\psi$ is onto.

Suppose $([v_1],[v_2])\in F^1$. Then \eqref{eqn_Soverlap} and
$S^+\cap(S^++e_1)=S\cap(S+e_1)$ imply
\[v_1^+(i)=v_2^+(i-e_1) \text{ for } i\in S^+\cap(S^++e_1).\]
Now Proposition~\ref{prop_skeleton} implies that there is a unique edge in
$B\Lambda(S^+,q,0,1)$ from $\psi^0([v_2])=v_2^+$ to
$\psi^0([v_1])=v_1^+$. Define a map $\psi^1$ from $F^1$ to the
edge set of $B\Lambda(S^+,q,0,1)$ by taking $\psi^1(([v_1],[v_2]))$ to be the
unique edge in $B\Lambda(S^+,q,0,1)$ with source $v_2^+$ and range
$v_1^+$. Then $\psi^1$ is a bijection since $\psi^0$ is. We have
$r\circ\psi^1=\psi^0\circ r$ since
\[
r(\psi([v_1],[v_2]))= \psi^0([v_1])= \psi^0(r([v_1],[v_2]))
\]
and similarly $s\circ\psi^1=\psi^0\circ s$. Thus $\psi=(\psi^0,\psi^1)$ is a graph isomorphism from $F$ to
$B\Lambda(S^+,q,0,1)$.
\end{proof}

So the blue graph of $T$ is related to the blue graph of $S^+$,
which is a tile with one fewer column than $T$. In fact we can
repeatedly apply Proposition~\ref{prop_dual} since the new tile
$S^+$ satisfies the hypotheses of that proposition. The tile $S_i^+$ in the next proposition is obtained from $T$
by deleting the first $i$ columns, shifting to
the origin and adding one box to the new first column.

\begin{corollary}\label{cor_tilesequence}
Suppose that $h_0>h_1$. For $1\leq
i\leq c_1$, let $S_i^+$ be the tile
\[S_i^+=((S^+_{i-1}\cap(S^+_{i-1}-e_1)\cup\{(h_i+1)e_2\},\]
define $r_{B_i}$ by
\[r_{B_i} =\begin{cases}
q^{h_0-h_1-1} &\text{ if } i=1\\
q^{h_{i-1}-h_i} &\text{ if } i>1,
\end{cases}\]
and let $H_i$ be the directed graph $B\Lambda(S_i^+,q,0,1)$. Then
$B\Lambda(T,q,t,w)\cong \widehat{(r_{B_1}H_1)}$ and $H_i\cong
\widehat{(r_{B_{i+1}}H_{i+1})}$ for $1\leq i\leq c_1-1$.
\end{corollary}

\begin{proof}Applying Proposition~\ref{prop_dual} to $B\Lambda(T,q,t,w)$
gives the result for $i=1$. Let $1\leq i \leq c_1-1$. Each tile $S_i^+$
has columns $h_{i}+1,h_{i+1},\ldots,h_{c_1}$. Since $T$ is hereditary, $h_i\geq h_{i+1}$. Then $S_i^+$ satisfies $(h_i+1)-h_{i+1}>0$
and so we can apply Proposition~\ref{prop_dual} to $H_i$ to get the result for $i>1$\end{proof}

\begin{examples}
(1) Suppose $T$ is the tile with $c_1=c_2=3$ and columns $h_0=3$,
$h_1=h_2=1$, $h_3=0$. Let $q=2$, $t\in\Zmod{2}$ and $w$ is a rule with
invertible corners. The tiles in Corollary~\ref{cor_tilesequence} are
\[T=
\beginpicture 
\setcoordinatesystem units <0.2cm,0.2cm> \setplotarea x from 0 to 4,
y from 0 to 4 \setlinear \plot 0 0 4 0 / \plot 0 0 0 4 / \plot 0 1 4
1 / \plot 0 2 3 2 / \plot 0 3 1 3 / \plot 0 4 1 4 / \plot 1 0 1 4 /
\plot 2 0 2 2 / \plot 3 0 3 2 / \plot 4 0 4 1 /
\endpicture \hspace{1cm}
S_1^+=\threestair \hspace{1cm} S_2^+=\bootup \hspace{1cm}
S_3^+=\twodominoup
\]
and the constants are $r_{B_1}=2$, $r_{B_2}=1$, $r_{B_3}=2$.

(2) Suppose $T$ has columns $h_0,\ldots,h_{c_1}$
satisfying $h_0=h_1+1$ and $h_1=h_2=\cdots=h_{c_1}$. Since
$S_{c_1}^+$ is the tile with one column with $h_{c_1}+1$ boxes, we
have $S_{c_1}\cap(S_{c_1}+e_1)=\emptyset$. So in $H_{c_1}$ there is
a directed edge between every pair of vertices, that is, $H_{c_1}$
is the complete graph $K_{q^{h_{c_1}}}$ with $q^{h_{c_1}}$ vertices.
Corollary~\ref{cor_tilesequence} implies that
$r_{B_1}=r_{B_2}=\cdots=r_{B_{c_1}}=1$ and so $B\Lambda(T,q,t,w)$ is
obtained by $c_1$ times taking the dual of $K_{q^{h_{c_1}}}$. For
example, the sock tile has $c_1=1$ and columns $h_0=1$ and $h_1=0$, and
so the blue graph of the Ledrappier graph is isomorphic to $\widehat{K_2}$.
\end{examples}

We need two more lemmas for the proof of
Proposition~\ref{prop_kernel0}.

\begin{lemma}\label{lemma_raeburn}
If $n\in\Z$ with $n>1$ and $B$ is an integer matrix, then
\[\ker(1-nB^t)=\{0\}.\]
\end{lemma}

\begin{proof}
Suppose $v\in\ker (1-nB)$, that is, $nBv=v$. We claim that $n^p|v$
for all $p\geq 1$. To see this, in the $p=1$ case we have $v=nBv$
and so $v$ is $n$ times some vector $Bv\in\Z^{\Lambda^0}$. For the
inductive step suppose $n^p|v$. Then there exists
$u\in\Z^{\Lambda^0}$ such that $v=n^pu$. Then
\[v=nBv=nB(n^pu)=n^{p+1}Bu\]
and so $n^{p+1}|v$. Hence $n^p|v$ for all $p\geq 1$, which is only
possible if $v=0$.
\end{proof}

\begin{lemma}\label{lemma_completekernel}
Let $n>1$ be an integer. If $K$ is the $n\times n$ matrix of all
$1$s, then \[\ker(1-K^t)=\{0\}.\]
\end{lemma}

\begin{proof}
The matrix $1-K^t$ is the circulant matrix\footnote{A circulant matrix is a square matrix in which each row vector is
obtained from the previous row vector by rotating one element to the
right. Then an $n\times n$ circulant matrix $\Circ(v)$ can be fully
specified by the first row vector
$v=(v_0,v_1,\ldots,v_{n-1})\in\Z^n$. (See \cite{Davis}).} $\Circ(v)$ with
$v=(0,1,\ldots,1)\in\Z^n$. If $\omega$ is a primitive $n$th root of
unity then using the formula for determinant of a circulant given in \cite{GR}
we have
\begin{align*}
\det (1-K^t) &= \det \Circ(v)
=\prod_{j=0}^{n-1}\sum_{i=0}^{n-1}\omega^{ij}v_i =
\prod_{j=0}^{n-1}\sum_{i=1}^{n-1}\omega^{ij}\\& =
\prod_{j=1}^{n-1}\sum_{i=1}^{n-1}\omega^{ij} \times
\sum_{i=1}^{n-1}\omega^{i0}
= \prod_{j=1}^{n-1}(-1)\times \sum_{i=1}^{n-1}1\\
&= (-1)^{n-1}(n-1).
\end{align*}
In particular $\det (1-K^t)\neq 0$, so we have $\ker(1-K^t)=\{0\}$.
\end{proof}

\begin{proof}[Proof of Proposition~\ref{prop_kernel0}]
We can deduce $(2)$ by applying part $(1)$ to the conjugate tile, so it suffices to prove $(1)$. Choose $r_{B_1},\ldots,r_{B_{c_1}}$ and
$B\Lambda,H_1,\ldots,H_{c_1}$ as in Corollary~\ref{cor_tilesequence}. Let $B,B_1,\ldots,B_{c_1}$ be the vertex matrices of
$B\Lambda,H_1,\ldots,H_{c_1}$. Since the vertex matrix of a dual graph $\widehat E$ is the edge matrix of $E$, Proposition~4.1 of \cite{MRS} gives isomorphisms
\begin{equation}\label{eqn_dualkerseq}
\ker(1-B^t) \cong \ker(1-r_{B_1}B_1^t) \text{ and } \ker(1-B_i^t)
\cong \ker(1-r_{B_{i+1}}B_{i+1}^t),
\end{equation}
for $1\leq i\leq c_1-1$. If $r_{B_i}=1$ for all $i$, then since the
tile $S_{c_1}^+$ has only one column,
Proposition~\ref{prop_skeleton}(b) implies that every entry in the
matrix $B_{c_1}$ is $1$, $\ker(1-B_{c_1}^t)=\{0\}$ by
Lemma~\ref{lemma_completekernel}, and all the kernels in
\eqref{eqn_dualkerseq} are trivial. If there exists $r_{B_j}$ which
is bigger than $1$, then there is a first such $j$; then
Lemma~\ref{lemma_raeburn} implies $\ker(1-r_{B_j}B_j^t)=\{0\}$, and
$\ker(1-B_i^t)=\{0\}$ for $i<j$. Hence $\ker(1-B^t)=\{0\}$.
\end{proof}

This completes the proof of Theorem~\ref{thm_kerdelta2trivial}, and
hence settles Conjecture~$(1)$. The next theorem settles
Conjecture~$(2)$.

\begin{theorem}\label{thm_K0eqK1}
Suppose we have basic data $(T,q,t,w)$ in which $w$ has invertible
corners and $c_1,c_2\geq 1$. Suppose further that  either
$h_0>h_1$ or $w_0>w_1$. Then the $\cstar$-algebra of the $2$-graph $\Lambda=\Lambda(T,q,t,w)$ has
\[|K_0(\cstarlambda)|=|K_1(\cstarlambda)|.\]
\end{theorem}

\begin{proof}
Since $h_0>h_1$ and $w_0>w_1$ we know from
Proposition~\ref{prop_kernel0} that $\ker(1-B^t)$ and $\ker(1-R^t)$
are trivial. Hence $1-B^t$ and $1-R^t$ are invertible over $\Q$, and both $K_0$ and $K_1$ are finite. Let
$C:=1-B^t$ and $D:=1-R^t$. Then $C$ and $D$ commute because $B$ and $R$ do, and
\[\delta_1=\begin{pmatrix}C &
D\end{pmatrix}:\ZL0\oplus\ZL0\to\ZL0 \text{ and } \delta_2=\begin{pmatrix} -D \\
C
\end{pmatrix}:\ZL0\to\ZL0\oplus\ZL0.\]
By Theorem~\ref{thm_kerdelta2trivial} we have $\ker\delta_2=\{0\}$
and $K_0(\cstarlambda)=\coker\delta_1$, that is,
\[K_0(\cstarlambda)=\ZL0/(C\ZL0+D\ZL0).\]
On the other hand we have \begin{align*}
K_1(\cstarlambda)&=\ker\delta_1 / \img\delta_2\\
&=\{(u,v):u,v\in\ZL0,Cu+Dv=0\} / \{(-Dw,Cw):w\in\ZL0\}.
\end{align*}
The map $C\ZL0\cap D\ZL0\to \ker\delta_1$ defined by $w\mapsto
(-C^{-1}w,D^{-1}w)$ carries $CD\ZL0$ onto
$\img\delta_2$. This induces an
isomorphism of $(C\ZL0\cap D\ZL0)/CD\ZL0$ onto
$\ker\delta_1/\img\delta_2$, hence
\[K_1(\cstarlambda)=(C\ZL0\cap D\ZL0) / CD\ZL0.\]
We have $C\ZL0\leq (C\ZL0+D\ZL0)\leq\ZL0$ and $C(D\ZL0)\leq
(C\ZL0\cap D\ZL0) \leq D\ZL0$. Since $D$ is an isomorphism of $\ZL0$ onto $D\ZL0$ which carries $C\ZL0$ onto $CD\ZL0$, we have  $|D\ZL0:CD\ZL0|=|\ZL0:C\ZL0|$.
Then
\begin{align}\label{eqn_cancel}
|D\ZL0:C\ZL0\cap D\ZL0|\, &|C\ZL0\cap D\ZL0:CD\ZL0| \\
&= |D\ZL0:CD \ZL0|\nonumber \\
&= |\ZL0:C\ZL0|\nonumber \\
&= |\ZL0:(C\ZL0+D\ZL0)|\, |(C\ZL0+D\ZL0):C\ZL0|.\nonumber
\end{align}
The inclusion of $D\ZL0$ in $C\ZL0+D\ZL0$ induces an isomorphism
of $D\ZL0/(C\ZL0\cap D\ZL0)$ onto $(C\ZL0+D\ZL0)/C\ZL0$, and hence
\begin{equation}|D\ZL0:C\ZL0\cap
D\ZL0|=|(C\ZL0+D\ZL0):C\ZL0|.\label{eqn_pink}\end{equation} Equation
\eqref{eqn_pink} allows us to cancel in \eqref{eqn_cancel} and
obtain
\[|C\ZL0\cap D\ZL0: CD\ZL0|=|\ZL0:(C\ZL0+D\ZL0)|,\]
which gives the result.
\end{proof}

\begin{remark}
Notice that our proof does not give an explicit isomorphism between
$C\ZL0\cap D\ZL0/CD\ZL0$ and $\ZL0/(C\ZL0+D\ZL0)$, so we cannot
deduce that $K_0\cong K_1$, only that they have the same number of
elements.
\end{remark}

\end{document}